\documentclass[11pt]{amsart}

\allowdisplaybreaks[1]

\usepackage{amssymb}
\usepackage{amsmath}
\usepackage{amscd}

\theoremstyle{plain}
\newtheorem{theorem}{Theorem}
\newtheorem{lemma}{Lemma}
\newtheorem{prop}{Proposition}
\newtheorem{coro}{Corollary}

\theoremstyle{definition}

\newtheorem{remark}{Remark}

\newcommand{\ts}{\hspace{0.5pt}}

\newcommand{\CC}{\mathbb{C}\ts}
\newcommand{\RR}{\mathbb{R}\ts}

\newcommand{\ZZ}{\mathbb{Z}}
\newcommand{\NN}{\mathbb{N}}

\newcommand{\TT}{\mathbb T}

\newcommand{\E}{\mathbb E}

\newcommand{\Rd}{{\mathbb R}^d}

\newcommand{\chF}{{\boldsymbol 1}}

\newcommand{\gl}{\lambda}
\newcommand{\gL}{\varLambda}

\newcommand{\CalA}{\mathcal{A}}
\newcommand{\CalB}{\mathcal{B}}
\newcommand{\CalL}{\mathcal{L}}

\newcommand{\CalE}{\mathcal{E}}

\newcommand{\CalM}{\mathcal{M}}
\newcommand{\CalO}{\mathcal{O}}
\newcommand{\CalS}{\mathcal{S}}
\newcommand{\CalX}{\mathcal{X}}
\newcommand{\Dm}{{\mathcal{D}}^{(m)}}
\newcommand{\cDr}{{\mathcal{D}}_r}

\newcommand{\X}{\ddot{X}}

\newcommand{\supp}{\mathrm{supp}}

\newcommand{\dmu}{{\mathrm d}\mu}
\newcommand{\card}{{\mathrm{card}}}
\newcommand{\iFT}{^\vee}
\newcommand{\norm}[1]{\left\Vert#1\right\Vert}
\newcommand{\abs}[1]{\left\vert#1\right\vert}
\newcommand{\Vol}{\ell}
\newcommand{\uF}{\boldsymbol{F}}

\newcommand{\FrP}{\mathfrak{P}}
\newcommand{\lfrq}{\underline{{\rm freq}}}
\newcommand{\ufrq}{\overline{{\rm freq}}}
\newcommand{\frq}{{\rm freq}}
\newcommand{\bm}{\boldsymbol{m}}
\newcommand{\palm}{\dot{\mu}}
\newcommand{\Av}{{\rm Av}}

\begin{document}

\title[Dworkin's argument revisited]
{
Dworkin's argument revisited:\\
Point Processes, Dynamics, Diffraction, \\
and Correlations
}

\author{Xinghua Deng and Robert V.\ Moody}
\address{Department of Mathematics and Statistics,
University of Victoria, \newline
\hspace*{12pt}Victoria, BC, V8W3P4, Canada}
\email{rmoody@uvic.ca, xdeng@math.ualberta.ca}
\date{\today}
\thanks{RVM gratefully acknowledges
the support of this research by the Natural Sciences and Engineering Research
Council of Canada.}

\begin{abstract}
  The setting is an ergodic dynamical system $(X,  \mu)$ whose points are themselves uniformly discrete point sets $\gL$ in some space $\Rd$ and whose group action is that of translation of these point sets by the vectors of $\Rd$. Steven Dworkin's argument relates the diffraction of the typical point sets comprising $X$ to the dynamical spectrum of $X$. In this paper we look more deeply at this relationship, particularly in the context of point processes.
  
We show that there is an $\Rd$-equivariant, isometric embedding, depending on the 
scattering strengths (weights) that are assigned to the points of $\gL \in X$, 
that takes the $L^2$-space of $\Rd$ under the diffraction measure into $L^2(X,\mu)$. We examine the image of this embedding and give a number of examples that show how it fails to be surjective. We show that
full information on the measure $\mu$ is available from the weights and set of {\em all} the correlations (that is, the $2$-point, $3$-point, \dots, correlations) of the typical point set $\gL \in X$.  

We develop a formalism in the setting of random point measures that includes multi-colour point sets, and arbitrary real-valued weightings for the scattering from the different colour types of  points, in the context of Palm measures and weighted versions of them. As an application we give a simple proof of a square-mean version of the  Bombieri-Taylor conjecture, and from that  we obtain an inequality that gives a quantitative relationship between the autocorrelation, the diffraction, and the $\epsilon$-dual characters of typical element of $X$. 
The paper ends with a discussion of the Palm measure in the context of defining pattern frequencies. \end{abstract}
\maketitle

\section{Introduction}
Imagine a point set representing the positions
of an infinite set of scatterers in some idealized solid of possibly infinite
extent. In practice such a set would be in $2$ or $3$ dimensional space, but
for our purposes we shall simply assume that it lies in some Euclidean space
$\Rd$. Suppose this point set satisfies the hard core condition that there
is a positive lower bound $r$ to the separation distance between the
individual scatterers (uniform discreteness). Consider the set $X$ of all possible
configurations $\gL$ of the scatterers. Assume that $X$ is invariant under the
translation action of $\Rd$ and assume also that there is a 
translation invariant ergodic probability measure $\mu$ on  $X$. 
In \cite{Dworkin}, Steven Dworkin pointed out
an important connection between the spectrum of the dynamical system
$(X, \Rd,\mu)$ and the diffraction of the scattering sets $\gL \in X$.

Dworkin's argument, as it is called (see Cor.~1 of Thm.~\ref{main}, below), has proven to be very fruitful, particularly in the
case of pure point dynamical systems and pure point diffraction,
where his argument for making the connection can be  made
rigourously effective, see for example  \cite{Hof,LMS1,Martin,BL}. Nonetheless,
the precise relationship between the diffraction and dynamics is quite elusive. 
One of the purposes of this paper is to clarify this connection.

The diffraction of
$\gL$, which is the Fourier transform of its autocorrelation (also
called the $2$-point correlation), is not
necessarily the same for all $\gL \in X$, whereas there is only one obvious measure, namely
$\mu$, on the dynamical system side with which to match it. However, the autocorrelation of $\gL$ is the same for $\mu$-almost 
all $\gL\in X$. In fact, as Jean-Baptiste Gou\'er\'e \cite{Gouere} has pointed out, using
concepts from  the theory of point processes there is a canonical construction
for this almost-everywhere-the-same autocorrelation through the use of the associated Palm measure $\palm$ of $\mu$.  Under the hypotheses above, the first moment $\palm_1$ of the Palm measure, which is a measure on the ambient space $\Rd$, is $\mu$-almost surely the autocorrelation of $\gL \in X$. We offer another proof of this in Theorems \ref{NautoCorr} and \ref{colourNautoCorr} below.

Put in these terms we can see that what underlies Dworkin's argument 
is a certain isometric embedding $\theta$ of the Hilbert space $L^2(\RR^d, \widehat {\palm_1})$ into $L^2(X, \mu)$. Both these Hilbert spaces afford natural unitary representations of $\RR^d$, call them $U_t$ and $T_t$ respectively ($t\in \RR^d$). Representation $T$ arises from the translation action of $\RR^d$ on $X$ and $U$ is a multiplication action which we define in $\eqref{Urep}$. The embedding $\theta$ intertwines the representations. However, $\theta$ is not in general surjective, and in fact it can fail to be surjective quite badly. 

The fact is that the diffraction, or equivalently the
autocorrelation measure of a typical point set $\gL \in X$, does 
not usually contain enough information to determine the measure $\mu$, even
qualitatively, see for example an explicit discussion of this in \cite{vEM}. We will give a number of other examples which show that outside the situation of pure point diffraction, one must assume that this is the normal state of affairs. In fact, even in the pure point case, $\theta$ can fail to be surjective. However, we shall show  in Thm.~\ref{correlationsDetermineMu} that, pure point or not, the knowledge of {\em all} the correlations of  $\gL$ ($2$-point, $3$-point, etc.) is  enough to determine $\mu$. This is one of the principal results of the paper and depends very much on the assumption of uniform discreteness. 

From a more realistic point of view, a material solid will be constituted from a number of different types of atoms and these will each have their own scattering strengths. We have incorporated this possibility into the paper by allowing there to be different types of points, labelled by indices $1,2, \dots, m$, and allowing each type (or colour, as we prefer to say)
to have its own scattering weight $w_i$. There is an important distinction to be made here. We view the measure $\mu$ on $X$ as depending only on the geometry of the point sets
(including the colour information) but not on the scattering weights, which only come into consideration of the diffraction. Thus $(X,\Rd,\mu)$ is independent of the weighting
scheme $w$, whereas the diffraction is not. The effect of this is that the diffraction is
described in terms of a weighted version of the first moment of the Palm measure
and the embedding $\theta^w: L^2(\RR^d, \widehat {\palm^w_1})\longrightarrow L^2(X, \mu)$
depends on $w$ (Thm.~\ref{main}). This allows us to study the significant effect that weighting has on this mapping.

As we have already suggested, an interesting and revealing point of view is to 
consider our dynamical system $(X,\Rd,\mu)$ as a point process, in which case
we think in terms of a random variable whose outcomes are the various point sets
$\gL$ of $X$. The theory of point processes is very well developed and has its own philosophy and culture. Although the theory is perfectly applicable to the situation that we are considering,
this  is nonetheless an unusual setting for it. More often random point processes are built around some sort of renewal or branching processes and the points sets involved do not satisfy any hard core property like the one we are imposing. Moreover,  diffraction is not a central issue. From the point of view of the theory of diffraction of tilings or Delone point sets, which are often derived in completely deterministic ways, it is not customary to think of these in terms of random point processes. But the randomness is not in the individual point sets themselves (though that is not disallowed, e.g. \cite{Hof4} or \cite{BBM}) but rather in the manner in which we choose them from $X$ and the way in which the measure $\mu$ of the dynamical system on $X$ can be viewed as a probability measure. The primary building blocks of the topology on $X$
are the cylinder sets $A$ of point sets $\gL$ that have a certain colour pattern in a certain finite region of space, and $\mu(A)$ is the probability that a point set $\gL$,
randomly chosen from $X$, will lie in $A$. 

The main purposes of this paper can be seen as continuing to build bridges between the study of uniformly discrete point sets (in the context of long-range order) and point processes that was started by Gou\'er\'e in \cite{Gouere}, and to provide a formalism of sufficient generality that the diffraction of point sets and the dynamics of their hulls can be studied together.

It is standard in the theory of point processes to model the point sets $\gL$ involved as
point measures $ \gl = \sum_{x\in \gL} \delta_x$, so that it is the supports of the
measures that correspond to the actual point sets. This turns out to be very convenient
for several reasons. The most natural topology for measures, the vague topology,
exactly matches the natural topology (local topology), which is used for the construction of
dynamical systems in the theory of tilings and Delone point sets (Prop.~\ref{vague=local}).  Ultimately, to discuss diffraction, one ends up in measures and the vague topology anyway, so having them from the outset is useful.  It is easy to build in the notion of colouring and weightings into measures. 

In fact there are a number of precedents for the study of
diffraction in the setting of measures rather than point sets \cite{BM,BL}. However, we note that the way in 
which weightings are used here does not allow one to simply start from arbitrary
weighted point measures at the outset. The point process itself knows about colours
but nothing about weights. As we have pointed out, the weighting only enters with the
correlations. 

The paper first lays out the basic notions of point processes, Palm measures, and moment measures, leading to the first embedding result, Prop.~\ref{link2}, mentioned above. We have chosen to develop this in the non-coloured version first, since this allows the essential ideas to be more transparent. The additional complications of colour and weightings are relatively easy to add in afterwards, leading to the main embedding theorem Thm.~\ref{main}, which establishes a mapping $\theta^w: L^2(\Rd,\widehat{\palm^w_1}) \rightarrow
L^2(X,\mu)$, $\palm^w_1$ being the weighted autocorrelation. The key to proving that the knowledge of all the higher correlations is enough to completely determine the law of the process is based on a result that shows that although $\theta^w$ need not be surjective, the algebra
generated by the image under $\theta^w$ of the space of rapidly decreasing functions on $\RR^d$  is
dense in $L^2(X,\mu)$ (with some restrictions on the weighting system $w$). This is
Thm.~\ref{AlgebraGenerates}. Here the uniform discreteness seems to play a crucial role.

Section~\ref{Examples} provides a number of examples that fit into the setting discussed here and that illustrate a 
variety of things that can happen. The reader may find it useful to consult this section in advance, as the paper proceeds.

As an application of our methods, we give a simple proof of a square-mean version of the  Bombieri-Taylor conjecture\footnote{For more on the history of this see Hof's discussion of it in \cite{Hof3}.}(see Thm.~\ref{Bombieri-Taylor}). Using this we obtain, under the assumption of finite local complexity, an inequality that gives a quantitative relationship between three fundamental notions: the autocorrelation, the diffraction, and the $\epsilon$-dual characters of typical elements of $X$. The proof of this does not involve colour and depends only on the embedding theorem Prop.~\ref{link2}. 

The paper ends with a discussion of the problem of defining pattern frequencies for elements of $X$, which arises because of the built-in laxness of the local
topology. We will find that the Palm measure provides a solution to the problem.

\section{Point sets and point processes}

\subsection{Point sets and measures} Start with $\Rd$ with its usual topology, and metric given by the Euclidean distance $|x-y|$ between points $x,y \in \Rd$. We let $B_R$ and $C_R$ denote the open ball of radius $R$ and the open cube of edge length $R$ about $0$ in $\Rd$. Lebesgue measure will be indicated by $\Vol$.

We are interested in closed discrete point sets $\gL$ in $\Rd$, but, as explained in the Introduction, we wish also to be able to deal with different types, or colours, of points. Thus we introduce 
$\bm := \{1, \dots, m\}$,
$m = 1,2,3, \dots$ with the discrete topology and take as our basic space the set $\E:= \Rd \times \bm$ with the product topology, so that any point $(x,i) \in \E$ refers to the point $x$ of $\Rd$ with colour $i$. When $m=1$ we simply identify $\E$ and $\Rd$.

Closures of sets in $\Rd$ and $\E$ are denoted by overline symbols. The overline also represents complex conjugation in this paper, but there is little risk of confusion.

There is the natural translation action of $\Rd$ on $\E$ given 
by
\[ T_t: (t, (x,i)) \mapsto t+ (x,i) := (t+x,i) \, .\]

Given $\gL \subset \E$, and $B\subset \Rd$, we define
\begin{eqnarray} \label{basicOps}
 B+ \gL &:= \bigcup_{b\in B} T_b\gL  &\subset \E \nonumber\\
 B \,\cap \,\gL &:= \{ (x,i) \in \gL\,:\, x\in B \} &\subset \E\\
 \gL^\downarrow &:= \bigcup_{(x,i) \in \gL} \{x\}  &\subset \Rd\, .\nonumber
 \end{eqnarray}
 $\gL^\downarrow$ is called the {\bf flattening} of $\gL$.
 
 Let $\CalO := \{(0,1), \dots, (0,m)\} \subset \E$. Then  
 $C_R^{(m)}:= C_R + \CalO$ is a `rainbow' cube that consists of the union of the cubes $(C_R,i)$, $i=1,\dots, m$. Its closure is $\overline{C_R^{(m)}}$.

Let $r >0$. A subset $\gL \subset \E$ is said to be $r$-{\bf uniformly discrete} if for all
$ a \in \Rd$, 
\begin{equation} \label{uniformdiscrete}
\card ((a + C_r)\cap \gL) \le 1 \, .
\end{equation}
In particular this implies that points of distinct colours cannot coincide.\footnote{It is more customary to define $r$-uniformly discreteness by using balls rather than cubes. This makes no intrinsic difference to the concept. However, in this paper, we find that the use of cubes makes certain ideas more transparent and some of our constructions less awkward.}  The family of all the $r$-uniformly discrete subsets of $\E$ will be denoted by 
$\cDr^{(m)}$. 

As we have pointed out, it is not just individual discrete point sets that we wish to discuss, but rather translation invariant families of such sets which collectively can be construed as dynamical systems. 

A very convenient way to deal with countable point subsets $\gL$ of $\E$ and families of them is to put them into the context of measures by replacing them by  pure point measures, where the atoms correspond to the points of the set(s) in question. To this end we introduce the following objects on any locally compact space $S$ :

\begin{itemize}
\item $\CalS$, the set of all Borel subsets of $S$;
\item $\CalB(S)$, the set of all relatively compact Borel subsets of $S$;
\item $BM_c(S)$, the space of all bounded measurable $\CC$-valued functions of compact support on $S$;
\item $C_c(S)$, the continuous $\CC$-valued functions with compact support on $S$.
If $S$ is known to be compact, we can write $C(S)$ instead.
\end{itemize}

Following Karr \cite{Karr} we let $M$ denote the set of all {\bf positive Radon measures} on $\E$, that is all positive regular Borel measures $\gl$ on $\E$ for which
$\gl(A) <\infty$ for all $A\in \CalB(\E)$. Equivalently, we may view these measures 
as linear functionals on the space $C_c(\E)$. We give $M$ the vague topology. This is the topology for which a sequence $\{\gl_n\} \in M$ converges to $\gl \in M$ 
if and only if $\{\gl_n(f) \} \rightarrow \gl(f)$ for
all $f \in C_c(\E)$. This topology has a number of useful characterizations, some of which we give below. 

Within $M$ we have the subset $M_p$ of {\bf point measures} $\gl$, those for which 
$\gl(A) \in \NN$ for all $A\in \CalB$. (Here $\NN$ is the set of natural numbers,
$\{0,1,2,\dots\}$.) These measures are always pure point measures in the sense that they are countable (possibly finite) sums of delta measures:

\[ \gl = \sum a_x \delta_x, \quad \quad x \in \E, a_x \in \NN \, .\]

Within $M_p$ we also have the set $M_s$ of {\bf simple point measures} $\gl$, those satisfying $\gl(\{x\}) \in \{0,1\}$, which are thus of the form
\[ \gl = \sum_{x\in \gL} \delta_x \]
where the support $\gL$ is a countable subset of $\E$. Evidently for these measures, for $x\in \E$,
\[\gl(\{x\}) >0 \Leftrightarrow \gl(\{x\}) = 1  \Leftrightarrow x \in \gL \, .\]

The Radon condition prevents the support of a point measure from having accumulation points in $\E$. Thus, the correspondence $\gl \longleftrightarrow \gL$ provides a bijection between $M_s$ and the closed discrete point sets of $\E$, i.e.
the discrete point sets with no accumulation points. 
This is the connection between point sets and measures that we wish to use.\footnote{Note that the point sets that we are considering here are simple in the sense that the multiplicity of each point in the set is just $1$. However, it is not precluded that the same point in $\Rd$ may occur more than once in such a point set, though necessarily it would have to occur with different colours. Very soon, however, we shall also preclude this.}
We note that the translation action of $\Rd$ on $\E$ produces an action of 
$\Rd$ on functions by $T_tf(x) = f(T_{-t} x)$, and on the spaces $M, M_p, M_s$ of
measures by $(T_t( \gl))(A)= \gl(-t+A), (T_t \gl)(f) = \gl(T_{-t}(f))$
for all $A\in \CalB(\E)$, and for all measurable functions $f$ on $\E$.

Here are some useful characterizations of the vague topology and some of its properties. These are cited in \cite{Karr}, Appendix A and appear with proofs in  \cite {DV-J}, Appendix A2.

\begin{prop} \label{vagueTop}
(The vague topology)
\begin{itemize}
\item[(i)] For $\{\gl_n\}, \gl \in M$ the following are equivalent:
 \begin{itemize}
  \item[(a)] $\{\gl_n(f) \} \rightarrow \gl(f)$ for all $f \in C_c(\E)$ (definition of vague convergence).
  \item[(b)] $\{\gl_n(f) \} \rightarrow \gl(f)$ for all $f \in BM_c(\E)$ for which 
the set of points of discontinuity of $f$ has $\gl$-measure $0$.
  \item[(c)]  $\{\gl_n(A) \} \rightarrow \gl(A)$ for all $A \in \CalB(\E)$ for which $\gl$ vanishes on the boundary of $A$, i.e. $\gl(\partial A)=0$.
 \end{itemize}
\item[(ii)] In the vague topology, $M$ is a complete separable metric space and $M_p$ is a closed subspace. 
\item[(iii)] A subspace $L$ of $M$ is relatively compact in the vague topology if and only if for all $A \in \CalB(\E)$, $\{\gl(A)\,:\, \gl \in L\}$ is bounded, which 
again happens if and only if for all $f \in C_c(\E)$,  $\{\gl(f)\,:\, \gl \in L\}$ is bounded. 
\end{itemize}
\end{prop}

Note that $M_s$ is not a closed subspace of $M_p$:  a sequence of measures in 
$M_s$ can converge to point measure with multiplicities.

\begin{prop} \label{BorelForM}
(The Borel sets of $M$)
The following $\sigma$-algebras are equal:
\begin{itemize}
\item[(i)] The $\sigma$-algebra $\CalM$ of Borel sets of $M$ under the vague topology.
\item[(ii)] The $\sigma$-algebra generated by requiring that all the mappings $\gl \mapsto \gl(f)$, $f\in C_c(\E)$ are measurable.
\item[(iii)] The $\sigma$-algebra generated by requiring that all the mappings $\gl \mapsto \gl(A)$, $A \in \CalB(\E)$ are measurable.
\item[(iv)] The $\sigma$-algebra generated by requiring that all the mappings $\gl \mapsto \gl(f)$, $f\in BM_c(\E)$ are measurable.
\end{itemize}
\end{prop}
 
A measure $\gl\in M$ is {\bf translation bounded} if for all bounded sets $K\in \CalB(\E)$,
$\{\gl(a+K) \,:\, a \in \Rd \}$ is bounded. In fact, a measure is translation bounded if this condition holds for a single set of the form $K = K_0 \times \bm$ where $K_0 \subset \Rd$ has a non-empty interior. For such a $K$ and for any positive constant $n$, we define the space $M_p(K,n)$ of translation bounded measures $\gl \in M_p$ for which 
\[\gl(a+K) \le n \] for all $a\in \Rd$. Evidently $M_p(K,n)$ is closed if $K$ is open, and by
Prop.~ \ref{vagueTop} it is relatively compact, hence compact. See also \cite{BL}, where this is proved in a more general setting.

If $r>0$ then $M_p(C_r^{(m)},1)$ is the set of point measures $\gl$ whose support $\gL$ satisfies the uniform discreteness condition (\ref{uniformdiscrete}). In particular, $M_p(C_r^{(m)},1) \subset M_s$ and is compact.

If $\gl \in M_s$ is a translation bounded measure on $\Rd$  we shall often write
expressions like $\sum_{x\in B} \gl(\{x\})$ where $B$ is some uncountable
set (like $\Rd$ itself). Such sums only have a countable number of terms
and so sum to a non-negative integer if $B$ is bounded, or possibly to $+\infty$
otherwise.

\subsection{Point processes} 
By definition, a {\bf point process} on $\E$ is a measurable mapping 
\[ \xi: (\Omega, \CalA, P) \longrightarrow  (M_p, {\CalM}_p) \]
from some probability space into $M_p$ with its $\sigma$-algebra of Borel sets
$\CalM \cap M_p$.
That is, it is a random point measure. Sometimes, when $m>1$, it is called
a ${\bf multivariate}$ point process. The {\bf law} of the point process is the  
probability measure which is the image  $\mu := \xi(P)$ of $P$. The point process is {\bf stationary} if $\mu$ is invariant under the translation action of 
$\Rd$ on $M_p$.

Thus from the stationary point process $\xi$ we arrive at a measure-theoretical dynamical system $(M_p, \Rd, \mu)$. Conversely, any such system may be interpreted as a stationary point process (by choosing $(\Omega, \CalA, P)$ to be $(M_p, \Rd, \mu)$).

There is no indication in the definition what the support of the law $\mu$ of the process might look like. In most cases of interest, this will be something, or be inside something, considerably smaller. In the sequel we shall assume that we have a point process $\xi: (\Omega, \CalA, P) \longrightarrow  (M_p, \CalM_p)$ that satisfies the following conditions:

\begin{itemize}
\item[(PPI)] the support of the measure $\mu = \xi(P)$ is a closed subset $X$ of
$M_p(C_r^{m},1)$ for some $r>0$.
\item[(PPII)] $\mu$ is stationary and has positive intensity (see below for
definition).
\item[(PPIII)] $\mu$ is ergodic.
\end{itemize}

These are examples of what are called translation bounded measure dynamical systems in \cite{BL}, although it should be noted that there the space of measures is not restricted to point measures, or even positive measures. 
 
Obviously under (PPI) and (PPII), $X$ is compact, and $(X, \Rd, \mu)$ is both a measure-theoretic and a topological dynamical system. 

Condition (PPI) implies that the point process is simple and so we may identify the measures of the point process as the actual (uniformly discrete) point sets in $\E$ that are their supports. Write $\X$ for the subset of $\Dm_r$ given by the supports of the measures of $X$. We call a point process satisfying (PPI) (and (PPII)) a {\bf uniformly discrete (stationary) point process}. We will make considerable use of these two ways of
looking at a point process -- either as being formed of point measures or of uniformly
discrete point sets. 

The ergodic hypothesis eventually becomes indispensable, but for our initial results it is not required. Usually we simplify the terminology and speak of a point process $\xi$ and assume implicitly the accompanying notation $(X, \Rd, \mu)$ and so on. We denote the family of all Borel subsets of $X$ by $\CalX$.
 
A key point is that the vague topology on the space $X$ of a uniformly discrete point process is precisely the topology most commonly used in the study of point set dynamical systems \cite{RW}. Sometimes this is called the {\bf local topology} since it implies a notion of closeness that depends
on the local configuration of points (as opposed to other topologies that depend only on the long-range average structure of the point set).

The local topology is most easily described as the uniform topology on $\Dm_r$ generated by the 
entourages
\begin{equation} \label{localTop}
U(C_R, \epsilon) := \\
\{(\gL,\gL') \in \Dm_r \, : \,   C_R \,\cap  \,\gL \subset C_\epsilon + \gL',
 \quad  C_R \,\cap  \,\gL' \subset C_\epsilon + \gL \} \,,
\end{equation}
where $R,\epsilon$ vary over the positive real numbers. 

Note that in (\ref{localTop}), $\gL$ and $\gL'$ are subsets of $\E$ and we are using the conventions  of (\ref{basicOps}). 
Intuitively two sets are close if on large cubes their points can be paired, 
taking colour into account, so that they are all within $\epsilon$-cubes of each other. It is easy to see that $\Dm_r$ is closed in this topology.

Given any $\gL' \in \Dm_r$ we define the open set
\[ U(C_R, \epsilon)[\gL'] := \{\gL \in \Dm_r \,:\, (\gL,\gL')\in U(C_R,\epsilon)\} \,.\]

 \begin{prop} \label{vague=local} {\rm (See also \cite{BL})} Let $\xi$ be a uniformly discrete point process. Then under the correspondence $\gl \leftrightarrow \gL$ between measures in $X$ and the point sets in $\X$, the vague and local topologies are the same.
 \end{prop}

{\bf Proof:} Let $r$ be the constant of the uniform discreteness. Let $\{\gl_n\}$ be a sequence of elements of $X$ for which the corresponding sequence $\{\gL_n\} \subset \X$ converges in the local topology to some point set $\gL \in \Dm_r$. Choose any positive function $f \in C_c(\E)$ and suppose that its support is in $C_R$, and choose
any $\epsilon >0$. Let 
$N_0 := 1 + \sup_{n \in \NN} \gl_n(C_R)$ and find $\eta >0$ so that $\eta <r$ and for all $x,y \in C_R$,
\[ |x-y| < \eta \Longrightarrow |f(x) -f(y)| < \epsilon/{N_0} \, .\] 
Let $\{x_1, \dots, x_N\} = C_R  \cap \gL \subset \E $. Then for all large $n$, $C_R \cap \gL_n \subset
\{C_\eta + x_1, \dots C_\eta + x_N\}$ with exactly one point in each of these cubes. Then
\[|\gl_n(f) - \gl(f) | = \left |\sum_{y\in C_R\cap \gL_n} f(y) - \sum_{x\in C_R\cap \gL} f(x)\right | \le N_0\epsilon/{N_0} = \epsilon \, . \] 
Thus $\{\gl_n(f)\}  \rightarrow \gl(f)$, and since $f\in C_c(\E)$ was arbitrary, 
$\{\gl_n\} \rightarrow \gl$ in $X$.

Now, going the other way, suppose that $\{\gl_n\} \rightarrow \gl$ in $X$. Let $R>0$ and let $C_R  \cap \gL = \{x_1, \dots, x_N\} $.  Choose any 
$0 <\epsilon < r$, small enough that for all $i\le N$, $C_\epsilon + x_i
\subset C_R^{(m)}$, and let 
\[f_\epsilon := \sum_{i=1}^N \chF_{C_\epsilon+x_i} \, .\] 
Then $\{\gl_n(f_\epsilon)\} \rightarrow \gl(f_\epsilon) = N = \gl(C_R) \leftarrow \{\gl_n(C_R)\}$, so for all
$n>>0$, $\gl_n(f_\epsilon) = N = \gl_n(C_R)$ (see Prop.~ \ref{vagueTop}). Since each cube $C_\epsilon+x_i$ can contain at most one point of any element of $\Dm_r$, then for all $n>>0$, and for all $i\le N$, there is a $y_i^{(n)} \in (C_\epsilon + x_i)\cap \gL_n$. This accounts for all the points of $C_R\,\cap \, \gL_n$. Thus $\gL_n \in U(C_R, \epsilon)[\gL]$. This proves that $\{\gL_n\} \rightarrow \gL$. \qed

\smallskip

\begin{remark} Let $\xi$ be a uniformly discrete point process. By Prop.~ \ref{vague=local}, the two topological spaces $X$ and $\X$ are homeomorphic. In particular, $\X$ is
compact in the local topology (a fact that can be seen directly from its definition). The the $\sigma$-algebras of their Borel sets $\CalX$ and $\ddot{\CalX}$ are isomorphic and we obtain a measure $\mu_{\X}$ on $\X$. Geometrically it is often easier to work in $\X$ rather than $X$, and we will frequently avail ourselves of the two different points of view. 
Notationally it is convenient to use the same symbols $X$ and $\mu$ for both and 
to use upper and lower case symbols to denote elements from $X$ according to whether we are treating them as sets or measures.
\end{remark}

\section{The moments and counting functions}

In this section we work in the one colour case $m=1$. Thus $\E = \Rd$.
We let $\xi: (\Omega, \CalA, P) \longrightarrow  (X, \CalX)$ be a uniformly discrete stationary point process on $\E$ with law $\mu$. We assume that
$X \subset M_p(C_r,1) \subset M_s$.

According to Prop.~\ref{BorelForM}, for each $A\in \CalB(\Rd)$ and 
for each $f \in BM_c(\Rd)$, the mappings 
\begin{eqnarray} \label{countingFunctions}
N_A: M_p(C_r,1) &\longrightarrow \ZZ, \quad &N_A(\gl) = \gl(A)\\
N_f: M_p(C_r,1) &\longrightarrow \CC, \quad &N_f(\gl) = \gl(f)
\end{eqnarray}
are measurable functions on $M_p(C_r,1)$, and by restriction, measurable functions
on $X$. The first of these simply counts the number
of points of the support of $\gl$ that lie in the set $A$, and $N_f$ is its natural extension from sets to functions. Whence the name {\bf counting functions}. They may also be considered as functions on $M_p$.
They may also be viewed as functions on the space $X$ viewed as the space of corresponding point sets.

 Thus, for example, in this notation we have for all $f \in BM_c(\Rd)$,
 \begin{eqnarray} \label{first moment}
 \int_X  \gl(f) d\mu(\gl) &=&  \int_X N_f(\gl) d\mu(\gl)
  =   \int_{X} \sum_{x\in \Rd} \gl(\{x\})f(x) d\mu(\gl) \\
 &=& \int_{X} N_f(\gL) d\mu(\gL) =   \int_{X} \sum_{x\in \gL} f(x) d\mu(\gL) 
 \, . \nonumber
 \end{eqnarray}

This is the first moment of the measure $\mu$, henceforth denoted $\mu_1$. More
generally, the $n$th {\bf moments}, 
$n = 1,2, \dots$ of a finite positive measure $\omega$ on $X$ are the 
unique measures on $(\Rd)^n$ 
defined by
\begin{eqnarray*}
 \omega_n(A_1 \times \dots \times A_n) &=&  
 \int_X \gl(A_1) \dots \gl(A_n)d\omega(\gl)\\
  &=&  \int_X N_{A_1} \dots N_{A_n} d\omega\, ,
 \end{eqnarray*}
where $A_1, \dots A_n$ run through all $\CalB(\Rd)$. Alternatively,
for all $f_1, \dots, f_n \in BM_c(\Rd)$,
 \[ \omega_n((f_1 , \dots , f_n)) = 
 \int_X N_{f_1} \dots N_{f_n} d\omega\, .\]
Since $\omega$ is a finite measure and the values of $\gl(f) =
N_f(\gl)$ are uniformly bounded for any $f\in BM_c(\Rd)$ as $\gl$ runs over 
$X$, these expressions define translation bounded measures on $(\Rd)^n$. 

If the measure $\omega$ is stationary  (invariant under the translation action
of $\Rd$) then the $n$th moment of $\omega$ is invariant under the
action of simultaneous translation of all $n$ variables.
Thus, if the point process $\xi$ is also stationary then the first moment 
of the law of $\xi$ is invariant, hence a multiple of Lebesgue measure:
\begin{equation} \label{first momentStationary}
\mu_1(A) = \int_X \gl(A) \dmu(\gl) = I\,\Vol(A) \,.
\end{equation}
This non-negative constant $I$, which is finite because of our assumption of uniform discreteness,  is the expectation for the number of points per unit volume of $\gl$  in $A$ and  is
called the {\bf intensity} of the point process. We shall always assume (see PPII) that the intensity is positive, i.e. not zero.
 
 The meaning of $N_f$ can be extended well beyond $BM_c(\Rd)$. To make this extension we introduce
 the usual spaces $L^p$-spaces $L^p(\Rd,\Vol)$, $L^p(X,\mu)$ together with their norms which we shall indicate by $|| \cdot ||_p$ in either case. In fact, we need these only for $p=1,2$. We shall also make use of the sup-norms $||\cdot||_\infty$.
 
 \begin{prop}\label{fToNf} The mapping {\rm (\ref{countingFunctions})} uniquely defines a continuous mapping (also called $N$)
\begin{eqnarray*} N: L^1(\Rd, \Vol) &\longrightarrow& L^1(X, \mu)\\
f &\mapsto& N_f
\end{eqnarray*}
satisfying $\norm{N_f}_1 \leq \sqrt{2}I \norm{f}_1$. Moreover, for all $f\in L^1(\Rd, \Vol)$,
\[ N_f(\gl) = \gl(f) \quad \mbox{ for $\mu$ almost surely all} \;\gl \in X \,.\]
 
\end{prop}

{\sc Proof:}  Let $A\subset \Rd$ be a bounded and measurable set, let $\chF_A$  
be the characteristic function of $A$ on $\Rd$, and define $N_{\chF_A}$ on $X$ by  $N_{\chF_A}(\gl) = \gl(\chF_A)
= \gl(A) = N_A(\gl)$, see (\ref{countingFunctions}). From (\ref{first momentStationary}),
$||N_{\chF_A}||_1 = \int_{X} N_{A}(\gl) d\mu(\gl)
= I \,\Vol(A) = I ||\chF_A||_1$. This shows that the result holds $N$ defined on these
basic functions.

For simple functions of the form
$f=\sum_{k=1}^{n}c_k\chF_{A_k}$, where the sets $A_k \subset \CalB(\RR^d)$ are 
mutually disjoint and the $c_k \ge0$, define
\[
N_f = \sum_{k=1}^{n}c_k N_{\chF_{A_k}}
= \sum_{k=1}^{n}c_k N_{A_k}\, .
\]
Then 
\[
\norm{N_f}_1= \sum_{k=1}^{n}c_k \norm{{N}_{\chF_{A_k}}}_1
= \sum_{k=1}^{n} c_k \,I  \Vol(A_k) = I\,\norm{f}_1\,,
\]
and $N_f(\gl) = \gl(f)$ for all $\gl\in X$.

The extension, first to arbitrary positive measurable functions and then to arbitrary real valued functions $f$ goes in the usual
measure theoretical way, and need not be reproduced here.

Finally we use the linearity to go to complex-valued integrable $f$.
If $f=f_{r}+\sqrt{-1} f_{i}$ is the splitting of $f$ into real and imaginary parts, then 
$N_f = N_{f_{r}}+\sqrt{-1}N_{f_{i}}$, so
\[
\norm{N_f}_1 \le I
(\norm{f_{r}}_1+\norm{f_{i}}_1)\\
=I \int_{\Rd}(\abs{f_r}+\abs{f_i})d\, \Vol.\\
\]
Using the inequality
$
\left(\abs{f_r}+\abs{f_i} \right)^{2}\leq
2(\abs{f_r}^{2}+\abs{f_i}^{2}),
$
we have
\[
\norm{N_f}_1\leq
\sqrt{2}I\int_{\Rd}\sqrt{\abs{f_r}^{2}+\abs{f_i}^{2}}d\, \Vol=\sqrt{2}I
\norm{f}_1.
\]

It is clear that if $f$ and $g$ differ on sets of measure $0$ then
likewise so do $N_f$ and $N_g$, so this establishes the existence of
the mapping.
 \qed

\begin{prop}\label{convInNf}
Let $f_n$, $n=1,2,3, \dots$ and
$f$ be measurable $\CC$-valued functions on $\Rd$ with supports all
contained within a fixed compact set $K$. Suppose that
$||f_n||_\infty, ||f||_\infty < M$ for some $M>0$ and $\{f_n\}
\rightarrow f$ in the $L^1$-norm on $\Rd$. Then $\{N_{f_n}\}
\rightarrow N_f$ in the $L^2$-norm on $X$.
\end{prop}

{\sc Proof:} Because of the uniform discreteness, 
$\gl(K)$ is uniformly bounded on $X$ by a constant $C(K) >0$.
Then for $g = f$ or $g= f_n$ for some $n$, $|N_g(\gl)| < MC(K)$.
\begin{eqnarray*}
|| N_f - N_{f_n}||_2^2 &=& \int_X \abs{N_f(\gl)  - N_{f_n}(\gl)}^2 d \mu(\gl)\\
&\le& \int_X (\abs{N_f(\gl)}  + \abs{N_{f_n}(\gl)}) \abs{N_f(\gl)  - N_{f_n}(\gl)}d \mu(\gl)\\
&\le& 2M C(K) \int_X \abs{N_f(\gl)  - N_{f_n}(\gl)}d \mu(\gl) \,,
\end{eqnarray*}
which, by Prop.~\ref{fToNf}, tends to $0$ as $n\to\infty$.
\qed

\section{Averages, the Palm measure and autocorrelation: 1-colour case}
\label{1Colour}

In this section we work in the one colour case $m=1$. Thus $\E = \Rd$.
We let $\xi: (\Omega, \CalA, P) \longrightarrow  (X, \CalX)$ be a uniformly discrete stationary point process on $\E$ with law $\mu$.

\subsection{The Palm measure}
The {\bf Campbell measure} is the measure 
$c'$ on $\Rd \times X$, defined by
\begin{equation}
c'(B\times D) = \int_D \gl(B) \,{\rm d}\mu(\gl)
= \int_X \sum_{x\in B} \gl(\{x\}) \chF_D(\gl) \,{\rm d}\mu(\gl)
\end{equation}
for all $B \times D \in \CalE \times \CalX$.

We note that $c'$ is invariant with respect to simultaneous translation of its two
variables. By introducing the measurable mapping 
\[\phi:\Rd \times X  \longrightarrow \Rd \times X:
\quad (x,\gl) \mapsto (x, T_{-x} \gl) \]
we obtain a twisted version $c$ of $c'$, also defined on $\Rd \times X$ : 
\begin{eqnarray*}
c(B\times D) &=&  \int_X \sum_{x\in B} \gl(\{x\}) \chF_D(T_{-x}\gl) \,{\rm d}\mu(\gl)\\
&=& \int_{X} \sum_{x\in B \cap \gL} (\chF_{D})(-x +\gL) \,{\rm d}\mu(\gL) \, ,
\end{eqnarray*}
and this
is invariant under translation of the first variable:
\begin{eqnarray*}
c((t+B)\times D)
&=& \int_X \sum_{x \in (t+B)}\gl(\{x\})\chF_D(T_{-x}\gl) \,{\rm d}\mu(\gl)\\
&=&\int_X \sum_{y\in B} T_{-t}\gl(\{y\}) \chF_D((T_{-y}T_{-t}\gl)  \,{\rm d}\mu(T_{-t}\gl)\\\\
 &=&c(B\times D)\, ,
\end{eqnarray*}
using the translation invariance of $\mu$.

Hence for $D$ fixed, $c$ is a multiple $\palm(D)\ell(B)$ of Lebesgue
measure and $D \mapsto \palm(D) := c(B\times D)/\ell(B)$ is a
non-negative measure on $X$ that is independent of the choice of $B\in \CalB(\Rd)$
(assuming that $B$ has positive measure). This measure is called the {\em Palm measure}
of the point process. See \cite{DV-J} for more details.

\begin{eqnarray} \label{Palm}
\palm(D) &= \frac{1}{\Vol(B)}\int_X \sum_{x\in B} \gl(\{x\}) \chF_{D}(T_{-x} \gl) \,{\rm d}\mu(\gl)\\
&=\frac{1}{\Vol(B)}\int_{X} \sum_{x \in B \cap \gL }\chF_{D}(-x + \gL) \dmu(\gL) \, .
\nonumber
\end{eqnarray}

We note that $\palm (X) = \int_X \gl(B) d \mu(\gl) / \Vol(B) = I$, which is the intensity of the point process. Some authors normalize the Palm measure by the intensity in order to render it a probability measure, and then call this probability measure the Palm measure. We shall not do this. However, we note that
the normalized Palm measure is often viewed as being the conditional probability
\[ \frac{1}{I}\palm(D) = \mu(\{\gl \in D\}\,|\, \gl(\{0\}) = 1  \}) \, , \]
that is, the probability conditioned by the assumption that $0$ is in the support of the point measures that we are considering.  In fact the conditional probability defined in this way is meaningless in general since the probability that $\gl(\{0\}) \ne 0$ is usually $0$. But the intuition of what is desired is contained in the definition. Taking $B$ as an arbitrarily small neighbourhood of $0$, we see that in effect we are only looking at points of $\gl$ very close to $0$ and then translating $\gl$ so that $0$ is in the support. The result is averaged over the volume of $B$. 

If the point process falls into the subspace $X$ of $\CalM$ then the support of the Palm measure is also in $X$. However, the Palm measure is not stationary in general, since the translation invariance of $\mu$ has, in effect, been taken out. 

Its first moment, sometimes called the intensity of the Palm measure, is
\begin{eqnarray} \label{palmIntensity}
\palm_1: \palm_1(A) &=& \int_X \gl(A) d\palm(\gl) \quad \mbox{or equivalently}\quad\\
\palm_1(f) &=& \int_X \gl(f) d\palm(\gl)  = \int_X N_f(\gl) d\palm(\gl) \nonumber\, .
\end{eqnarray}

The first moment of the Palm measure, and also the higher moments to be defined later, play a crucial role in the development of the paper, since they are, in an almost sure sense, the $2$-point and higher point correlations of the elements of $X$.

As with $\mu$, we will, consider the Palm measure interchangeably as a
measure on $X$ or on  $\X$ (as we have already done implicitly in  
Eq. (\ref{Palm})). 

The importance of the Palm measure is its relation to the average value of a function over a typical point set $\gL \in X$, and from there to pattern frequencies in $\gL$ and its direct involvement in the autocorrelation of $\gL$. To explain this we need to develop the Palm theory a little further.

\begin{lemma} (Campbell formula) \label{Ncampbell}  For any measurable
function $F:\Rd \times X \longrightarrow \RR$,
\[\int_{\Rd} \int_{{X}} F(x, \gl) \textrm{d} \palm(\gl) \, {\rm d}x
= \int_{{X}} \sum_{x \in \Rd} \gl(\{x\})F(x, T_{-x}\gl) \, {\rm d}\mu (\gl) \, .\]
\end{lemma}

{\sc Proof:} This can be proven easily by checking it on simple
functions. Let $F=\chF_B\times\chF_D$. Then
\begin{eqnarray*}
\int_{\Rd}\int_{{X}}\chF_B(x)\times\chF_D(\gl)\textrm{d}
\palm(\gl){\rm d}x &=& \ell(B)\palm(D)=c(B\times D)\\
&=&\int_{{X}}\sum_{x\in \Rd} \gl(\{x\})\chF_B(x)\chF_D(T_{-x}\gl)d\mu(\gl)\\
&=& \int_{{X}} \sum_{x \in \Rd} \gl(\{x\})F(x, T_{-x} \gl) \, \dmu(\gl) \, .
\end{eqnarray*}
\qed

Let $\nu_R$ be the function on $X$ defined by
\[\nu_R(\gl)=\frac{1}{\Vol (C_R)}N_{C_R}(\gl)\, ,\]
for all $R>0$.
We treat $\nu_R$ as the Radon-Nikodym density of an absolutely continuous measure on $X$ (with respect to $\mu$).

\begin{prop}\label{limitmeasure}

In vague convergence, 
\[ \{\nu_R\}\rightarrow \palm \quad \mbox{as} \quad R \to 0  \, .\]
\end{prop}

{\sc Proof:} Use the definition of the Palm measure in (\ref{Palm}) with $B$
replaced by $C_R$. Then for any continuous function $G$ on $X$,
\[
\palm(G)=\frac{1}{\Vol(C_R)}\int_{X}\sum_{y\in C_R}\gl(\{y\})G(T_{-y}\gl)d\mu(\gl).
\]

If we require that $R<r$ then
\[\sum_{y\in C_R}\gl(\{y\})G(T_{-y}\gl)=N_{C_R}(\gl)G(T_{-x}\gl) \, ,\] where $x$ is the
unique point in $\gL \cap C_R$ when it is  not empty, and then
\[
\palm(G)=\frac{1}{\Vol(C_R)}\int_{X}N_{C_R}(\gl)G(T_{-x}\gl)d\mu(\gl).
\]

On the other hand
\[
\nu_R(G)=\frac{1}{\Vol(C_R)}\int_{X}N_{C_R}(\gl)G(\gl)d\mu(\gl).
\]

Thus
\begin{eqnarray}\label{difference}
&&\abs{\palm(G)-\nu_R(G)} \nonumber\\
&=&\abs{\frac{1}{\Vol(C_R)}\int_{X}N_{C_R}(\gl)\{G(T{-x}+\gl)-G(\gl)\}d\mu(\gl)}\\
&\leq&\frac{1}{\Vol(C_R)}\int_{X}N_{C_R}(\gl)\abs{\{G(T_{-x}\gl)-G(\gl)\}}d\mu(\gl)\,.
\nonumber
\end{eqnarray}

The rest follows from the uniform continuity of $G$ ($X$ is compact). From the inequality (\ref{difference}), 
\[
\abs{\palm(G)-\nu_R(G)}\leq
\frac{\epsilon_R}{\Vol(C_R)}\int_{X}N_{C_R}(\gL)d\mu(\gL)=\epsilon_R
I \, ,\
\]
for some $\epsilon_R \to 0$ as $R\to 0$, where $I$ is the intensity of the point process.

Therefore, we have that $\nu_R\rightarrow \palm$ vaguely. \qed

\subsection{Averages}

Let $\xi$ be a uniformly discrete ergodic stationary point process, with corresponding dynamical system $(X,\Rd, \mu)$. Let $F\in C(X)$. 
The {\bf average} of $F$ at $\gl\in X$ is  
\[\Av(F)(\gl) = \lim_{R\to\infty}\frac{1}{\Vol(C_R)} \sum_{x\in C_R}\gl(\{x\}) F(T_{-x}\gl) \, ,\]
if it exists. Thus $\Av(F)$ is a function defined at certain points of $X$.
Alternatively, we may think of $F$ as a function on point sets and write this as 

\[\Av(F)(\gL) = \lim_{R\to\infty}\frac{1}{\Vol(C_R)} \sum_{x \in \gL\cap C_R} F(-x + \gL) \, .\]
We will prove the almost-sure existence of averages.

\begin{prop} \label{averageThm}
Let $F\in C(X)$. The average value of $\Av(F)(\gl)$ of $F$  exists $\mu$-almost surely for
$\gl \in X$ and it is almost surely equal to $\palm(F)$. In particular $\Av(F)$ exists as a measurable function on $X$. 
If $\mu$ is uniquely ergodic then the average value always exists everywhere and is equal to
$\palm(F)$.
\end{prop}

{\sc \bf Proof:} It is clear that the average value is constant along the orbit of
any point $\gl$ for which it exists.

Let $\epsilon >0$. Since $F$ is uniformly continuous,
 there is a compact set $K$ and an $s>0$ so that 
$|F(\gl') - F(\gl'')| <\epsilon$ whenever $(\gl',\gl'')\in U(K,s)$. In particular
$|F(-x+\gl) - F(-u+\gl)|<\epsilon$ whenever $|x-u|<s$. We can assume that 
$s<r$. 

Let $\nu_s := \frac{1}{\Vol(C_s)}N_{C_s}  : X \longrightarrow \CC$, as above.
For $x \in \Rd$ and $\gl \in X$, 
$N_{C_s}(T_{-x}\gl) =1$ if and only if $x\in u + C_s$ for
some $u\in \gL$. Thus
\begin{align}
 \frac{1}{\Vol(C_R)} \int_{C_R} & F(T_{-x}\gl) \nu_s(T_{-x}\gl) dx \\
 &\sim 
\frac{1}{\Vol(C_R)} \sum_{u \in C_R} \gl(\{u\}) \frac{1}{\Vol(C_s)}\int_{u + C_s} F(T_{-x}\gl) dx \,,
\end{align}
where the $\sim$ comes from boundary effects only and becomes equality in the limit.

There is a constant $a>0$ so that $\card (\gl(C_R))/\Vol(C_R) < a$, independent of $R$
or which $\gl \in X$ is taken. Using this and our choice of $s$, we obtain

\begin{align}  \label{estimate}
\vert\lim_{R\to\infty}\frac{1}{\Vol(C_R)} \int_{C_R} &F(T_{-x}\gl) \nu_s(T_{-x}\gl) dx\\ 
 &- \lim_{R\to\infty} \frac{1}{\Vol(C_R)} \sum_{u \in  C_R}  \gl(\{u\}) F(T_{-u}\gl)  \vert < a\epsilon \,. \nonumber
\end{align}

The right hand term is just the average value of $F$ at $\gl$, if the limit exists. However, by the Birkhoff ergodic theorem the left integral exists almost surely and is equal to 
$\int_{X} F \nu_s \dmu = \nu_s(F)$, $\nu_s$ being treated as a measure.

Now making $\epsilon \to 0$, so $s\to 0$ also, and using Prop. \ref{limitmeasure} we have
\[ \palm(F) = \lim_{s\to 0} \nu_s(F) = \lim_{R\to\infty}\frac{1}{\Vol(C_R)} \sum_{u \in C_R} 
\gl(\{u\})F(T_{-u} \gl)  = \Av(F)(\gl)\, .\]
Thus the average value of $F$ on $\gl$ exists almost surely.

In the uniquely ergodic case, the conclusion of Birkhoff's theorem is true everywhere
in $X$.
\qed

\subsection{The autocorrelation and the Palm measure}

Again, let \newline
$\xi: (\Omega, \CalA, P) \longrightarrow  (X, \CalX)$ be a uniformly discrete stationary ergodic point process on $\Rd$ with law $\mu$. For each $\gl \in X$ we define $\tilde{\gl}$ to be the 
point measure on $\Rd$ defined by $\tilde{\gl}(\{x\}) = \overline{\gl(\{-x\})}$ (though
at this point we are only dealing with real measures). Then the autocorrelation of $\gl$ is defined as 

\[\gamma_{\gl} := \lim_{R\to \infty} \frac{1}{\Vol(C_R)} \left( \gl|_{C_R} \ast
\widetilde{\gl}|_{C_R} \right)= \lim_{R\to \infty}\frac{1}{\Vol(C_R)}\sum_{x,y\in \gL\cap
C_R}\delta_{y-x} \, .\]
Here the limit, which may or may not exist, is taken in the vague topology.

A simple consequence of the van Hove
property of cubes is:

\begin{equation}\label{AL1}
\gamma_{\gl} =
\lim_{R\to\infty}  \frac{1}{\Vol(C_R)}\sum_{x \in\gL\cap C_R,y\in\gL}\delta_{y-x} \, .
\end{equation}
Namely, for any $f\in C_c(\Rd)$, say with support $K$, and any $x\in C_R$, $f(y-x) =0$
unless $y \in C_R +K$, and thus for large $R$ the only relevant $y$ which are not
in $C_R$ are in the $K$-boundary of $C_R$, which is vanishingly small in relative volume as $R \to\infty$.

\begin{theorem}\label{NautoCorr}  The first moment $\palm_1$ of the Palm
measure is a positive, positive definite, translation bounded measure. Furthermore, $\mu$-almost surely, $\gl\in X$ admits an autocorrelation $\gamma_{\gl}$ and it is equal to $\palm_1$.
If $X$ is uniquely ergodic then  $\palm_1 = \gamma_{\gl}$ for all $\gl \in X$.
\end{theorem}

{\sc \bf Proof:}
We begin with the statement about the autocorrelation measures $\gamma_\lambda$. Let $f\in C_c(\Rd)$. The autocorrelation of $\gl$ at $f$, if it exists, is
\begin{align} \label{basicACequation}
\gamma_{\gl}(f) 
 &= \lim_{R\to\infty} \frac{1}{\Vol(C_R)}\sum_{x \in C_R, y\in \gL}\gl(\{x\}) f(y-x)\\\nonumber
 &=  \lim_{R\to\infty} \frac{1}{\Vol(C_R)} \sum_{x \in C_R} \gl(\{x\})N_f(T_{-x}\gl)\\\nonumber
 &= \palm(N_f) = \palm_1(f) \nonumber
\end{align}
for $\gl \in X$, $\mu$-almost surely,
where we have used Prop. \ref{averageThm} and (\ref{palmIntensity}).

This is basically what we want, but we must show that it holds for {\em all} $f\in C_c(\Rd, \RR)$
for almost all $\gl \in {X}$. This is accomplished by using a countable dense (in the
sup norm) set of elements of $C_c(\Rd, \RR)$. We can get (\ref{basicACequation}) simultaneously for this countable set, and this is enough to get it
for all $f\in C_c(\Rd, \RR)$. Then $\gamma_{\gl}$ exists and is equal to $\palm_1$ for almost
all $\gl \in {X}$. For more details see \cite{Gouere}. 

Finally, it is clear that $\gamma_\gl$ is a positive and positive definite measure whenever it exists,
and hence also $\palm_1$ is positive and positive definite. All positive and positive definite
measures are translation bounded, \cite{BF} Prop.~4.4., or \cite{Hof}.\qed

\subsection{Diffraction and the embedding theorem}

As a consequence of the positive and positive-definiteness of the autocorrelation 
 this that they are Fourier transformable and that their Fourier transforms
are likewise positive, positive definite, and translation bounded, \cite{BF}. 

We recall that the Fourier transform of such a measure $\omega$ on $\Rd$ can be 
defined by the formula:
\begin{equation}
\widehat{\omega}(f) = \omega(\widehat{f})
\end{equation}
for all $f$ in the space $\mathbb{S}$ of rapidly decreasing functions of $\Rd$. In fact, it will suffice to have this formula on the space $\mathbb{S}_c$ of compactly supported functions in $\mathbb{S}$, since they are dense in $\mathbb{S}$ in the standard topology on
$\mathbb{S}$ (\cite{Schwartz}). The key point is that if $\{f_n\} \in \mathbb{S}_c$ converges to
$f \in \mathbb{S}$, then $\{\widehat{f_n}\}$ converges to $\widehat{f}$ and one can use the translation boundedness of $\omega$ to see then that 
$\{\omega(\widehat{f_n})\}$ converges to $\{\omega(\widehat{f})\}$, i.e.
$\widehat{\omega}(f)$ is known from the values of $\{\widehat{\omega}(f_n)\}$.

The measure
$\widehat{\gamma_\gl}$ is the {\bf diffraction} of $\gl$, when it exists. Our results show that the first moment of the Palm measure, $\palm_1$ must also be a positive, positive definite 
transformable translation bounded measure and that almost surely $\widehat{\palm_1}$ is the diffraction of $\gl \in X$. 

The next result appears, in a slightly different form in \cite{Gouere}.
For complex-valued functions $h$ on $\E$ define $\tilde h$ by $\tilde h (x) = 
\overline{h(-x)}$. We denote the standard inner product defined by $|| \cdot ||_2$
on $L^2(X,\mu)$ by $(\cdot, \cdot)$.

\begin{prop} \label{link1}
Let $g,h\in BM_c(\RR^d)$ and suppose that 
$g*\tilde h *\palm_1$ is a continuous function on $\RR^d$.
Then for all $t\in \RR^d$,
\[g*\tilde h*\palm_1 (-t) = (T_t N_g, N_h) \, . \]
\end{prop}

{\sc Proof:} It suffices to prove the result when $g,h$ are real-valued functions.
By Prop.\ref{BorelForM},
$N_g, N_h$ are measurable functions on ${X}$, and they are clearly $L^1$-functions (Prop.~\ref{fToNf}).
\begin{eqnarray*}
g*\tilde h*\palm_1(-t) &=& \int_{\RR^d} (g*\tilde h)(-t-u) {\mathrm d}\palm_1(u) =
 \int_{\RR^d} (\tilde g*h)(t+u) {\mathrm d}\palm_1(u) \\
 &=&  \int_{\RR^d} \widetilde{T_t g}*h(u) {\mathrm d}\palm_1(u) \\
 &=&
 \int_{X} \left(\sum_{x \in \Rd} \gl(\{x\})(\widetilde{T_{t} g} * h)(x) \right)
{\mathrm d}\palm(\gl)\\
 &=& \int_{\RR^d}\int_{X} \sum_{x \in \Rd} \gl(\{x\})(T_{t}g)(u)h(x+u)  
\,{\mathrm d}\palm(\gl) {\mathrm d}u\\
 &=& \int_{\RR^d} \int_{{X} }(T_{t}g)(u)T_{-u} N_h(\gl) {\mathrm
d}\palm(\gl) \, {\mathrm d}u
  \end{eqnarray*}
where we have used (\ref{first moment}) and the dominated convergence theorem to
rearrange the sum and the integral.
  Now using the Campbell formula we may continue:
  \begin{eqnarray*}
\palm_1*g*\tilde h(-t) &=&
 \int_{{X} } \sum_{u\in \Rd} \gl(\{u\}) (T_{t}g)(u)N_h(\gl) \dmu(\gl)\\
&=& \int_{{X}}  N_{T_{t}g}(\gl) N_h(\gl) \dmu = (T_{t}N_g, N_h) \, .
 \end{eqnarray*}  \qed
 
We are now at the point where we can prove the embedding theorem (in the unweighted 
case). This involves the two Hilbert spaces $L^2(\Rd, \widehat{\palm_1})$ and $L^2(X, \mu)$. Since the translation action of $\Rd$ on $X$ is measure
preserving, it gives rise to a unitary representation $T$ of $\Rd$
on $L^2(X, \mu)$ by the usual  translation action of $\Rd$ on measures.

We also have a unitary representation $U$ of $\Rd$ on
$L^2(\Rd, \widehat{\palm_1})$ defined by
\begin{equation} \label{Urep}
U_tf(x) = e^{-2\pi i t.x} f(x) = \chi_{-t}(x)f(x)\, ,
\end{equation}
where the characters $\chi_k$ are defined by 
\begin{equation} \label{chars}
\chi_k(x) = e^{2 \pi i k.x}\, .
\end{equation}
We denote the inner product of $L^2(\Rd, \widehat{\palm_1})$
by $\langle \cdot, \cdot \rangle$ and note that with respect to it $U$
is a unitary representation of $\RR^d$.

\begin{prop}\label{link2}
If $g,h \in \mathbb{S}$  are rapidly decreasing functions then
\[g*\tilde h *\palm_1(-t)  = \langle U_t (\hat g), \hat h\rangle \, .\]
In particular,
\[ \langle U_t (\hat g), \hat h\rangle = (T_tN_g,N_h) \, . \]
Thus there is an isometric embedding intertwining $U$ and $T$, 
\[\theta :L^2(\RR^d, \widehat{\palm_1}) \longrightarrow L^2(X, \mu) \, ,\]
under which 
\[ \hat f \mapsto N_f \]
for all $f \in \mathbb{S}$.

\end{prop}

{\sc Proof:} As we have pointed out, it will suffice to show the first result for 
$g,h \in \mathbb{S}_c$ since it is dense in $\mathbb{S}$ under the standard
topology of $\mathbb{S}$. 
We note that the hypotheses of Prop.~\ref{link1} are satisfied, so, starting as in its proof and denoting the inverse Fourier
transform by $f \mapsto \check f$, we have
\[
g*\tilde h*\palm_1(-t) =\int_{\E} \widetilde{T_t g}*h(u) {\mathrm d}\palm_1(u)
= \int_{\E} (\widetilde{T_t g})\iFT h\iFT {\mathrm d}\widehat{\palm_1} \, .
\]
The first result follows from $\check{h} = \check{{\overline{h}}} =
\overline{\hat h}$ and $(\widetilde{T_t g})\iFT =
\widehat{T_t g} = \chi_{-t} \hat g$.  

The second part of the proposition follows from Prop.~\ref{link1} and the observation that  $\mathbb{S}_c$ is dense in $C_c(\Rd)$ in the sup norm (\cite{Schwartz}, Thm.~1), hence
certainly in the $||\cdot ||_2$-norm, and $C_c(\Rd)$ is dense in
$L^2(\Rd, \widehat{\palm_1})$ in the $||\cdot ||_2$-norm (see \cite{Rud}, Appendix E). 

Thus we have the existence of the embedding on a dense subset of 
$L^2(\RR^d, \widehat{\palm_1})$ and it extends uniquely to the closure.\qed

\section{Adding colour}

We now look at the changes required to Section \ref{1Colour} in order to include
colour, i.e. to have $m>1$. The colour enters in two ways. First of all, the dynamics,
that is to say the dynamical hull $X$ and the measure $\mu$, depend on colour since closeness in the local topology depends on simultaneous closeness of points of  like-colours. Secondly the autocorrelation, and then the diffraction, depends on colour. 

Diffraction depends on how scattering waves from different points (atoms) superimpose upon each other. However, physically, different types of atoms will have different scattering strengths, and so we wish to incorporate this into the formalism. This is accomplished by
specifying a vector $w$ of weights to be associated with the different colours and
introducing for each point measure $\gl$ of our hull $X$ a weighted version of
it, $\gl^w$. This will be a measure on $\Rd$. It will be important that the weighting is kept totally separate from the topology and geometry of $X$. The geometry of the configuration and the weighting of points, which enters into the diffraction, are different things. The measures describing our point sets are measures on $\E$, but the diffraction always takes
place on the flattened point sets. 

On the geometrical side we have treated the full colour situation from the start. In this section we introduce it into the autocorrelation/diffraction side. This affects almost every result in Section \ref{1Colour}. However, we shall see that every proof then generalizes quite easily, and we simply outline the new situation and 
the generalized results, leaving the reader to do the easy modifications to the proofs.

\subsection{Weighting systems}

Let $\xi: (\Omega, \CalA, P) \longrightarrow (X,\CalX)$ be a uniformly discrete
stationary multi-variate point process, where 
$X \subset M_p$, $\CalX = \CalM \cap X$, and $(X,\RR^d, \mu)$ is the resulting
dynamical system. We let $\E = \RR^d \times \bm$, with $\E^i = \RR^d \times \{i\}$
and $\E = \bigcup_{i\le m}\E^i$.
For each $i$ we have the restriction 
\[{\rm res}^i : \gl \mapsto \gl^i \]
of measures on $\E$ to measures on $\E^i$. We will simply treat these restricted measures as being measures on $\RR^d$. If $\gl \leftrightarrow \gL$ then we also
think of ${\rm res}^i$ as the mapping $\gL \mapsto \gL^i:= \{x \in \RR^d: (x,i) \in \gL\}$.
\footnote{It is also possible to define associated dynamical systems $X^i$ and with them Palm measures. However, it is important here that everything will always refer back to the full colour situation encoded in the geometry of $X$.}

The same argument that led to (\ref{first momentStationary}) gives
\begin{equation}\label{colourIntensity}
\int_X \lambda^i(A) d\mu(\lambda) = I^{(i)} \Vol(A) 
\end{equation}
for some $ I^{(i)} \ge 0$, for each $i$. We shall always assume:

\smallskip
\smallskip \noindent
\begin{itemize}
\item[(PPIIw)]  \quad $I^{(i)} >0$ for all $i\le m$.
\end{itemize}

\smallskip

A {\bf system of weights} is a vector $w = (w_1, \dots w_m)$ of real numbers.\footnote{One could have complex numbers here, but it makes things easier, 
and more natural for higher correlations, if the weights are real.} We define a mapping 
\[X \to M_s(\Rd)  \quad\gl \mapsto \gl^w := \sum_{i\le m} w_i \gl^i \, .\]
The quantity
\begin{equation}\label{weightedIntensity}
I^w := \sum_{i=1}^m w_i I^{(i)}
\end{equation}
is called the {\bf weighted intensity} of the weighted point process.

 We also have the flattening map:
 \[X \to M_s(\Rd)  \quad\gl \mapsto \gl^\downarrow := \sum_{i\le m} \gl^i \, .\]

First introduce the measure $c^w$ on  $\RR^d \times X$:
\[ c^w(B\times D) = \int_X\sum_{x\in B} \gl^w(\{x\})T_x \chF_D(\gl) {\rm d}\mu(\gl) \, .\]
Since $(T_x\gl)^w = T_x(\gl^w)$ this measure is invariant under translation of the first variable and we have 
\[ c^w(B\times D) = \Vol(B) \palm^w(D) \,.\]
This determines the $w$-{\bf weighted Palm measure} $\palm^w$ on $X$. This is not a Palm measure in the normal sense of the word. However, it plays the same role
as the Palm measure in much of what follows. For example, there is a corresponding Campbell formula: 
\[ \int_{\RR^d} \int_{X} F(x,\gl) {\rm d}\palm^w(\gl) {\rm d}x 
= \int_{X} \sum_{x \in \Rd} \gl^w(\{x\}) F(x ,T_{-x} \gl) {\rm d} \mu(\gl) \]
for all measurable $F:\RR^d \times X \longrightarrow \CC$.

We note the formula for the weighted intensity:
\begin{eqnarray} \label{intensityFormula}
I^w l(A) &=& \int_X\gl^w d \mu(\gl) = \int_X \sum_{x \in A} \gl^w(\{x\}) d \mu(\gl) \nonumber\\
&=& c^w(A \times X) = l(A) \palm^w(X) \, , \mbox{whence} \nonumber\\
I^w &=& \palm^w(X) \,.
\end{eqnarray}

For all $i\le m$, for all $A\in \CalB(\Rd)$, and for all $f\in BM_c(\RR^d)$ define
\begin{eqnarray} \label{weightedN}
N_A^w : X \longrightarrow \NN & \quad N_A^w(\gl) = \gl^w(A)  = \sum_{x\in A}\gl^w(\{x\})\\
N_f^w : X \longrightarrow \NN & \quad N_f^w(\gl) = \gl^w(f)  = 
\sum_{x\in \Rd}\gl^w(\{x\}) f(x)\, .
\nonumber
\end{eqnarray}
Thus, for example,
\begin{equation}\label{useOfres}
N_A^w(\gl) = \sum w_i \gl^i(A) = \sum w_i N_A( {\rm res}^i(\gl))
 = \sum w_i N_A\circ {\rm res}^i (\gl) \, .
\end{equation}

Define 
\[\nu_R^w: X \longrightarrow \RR, \quad \nu_R^w(\gl) = \frac{1}{\Vol (C_R)} 
N_{C_R}^w(\gl)\]
or equivalently, $\nu_R^w(\gL) = \frac{1}{\Vol (C_R)} N_{C_R}^w(\gL)$.
In vague convergence, 
\[ \{\nu_R^w\}\rightarrow \palm^w \quad \mbox{as} \quad R \to 0  \, .\]

These auxiliary measures are used, as before, to prove the existence
of averages.
Let $F \in C(X)$. The $w$-{\bf {average value}} of $F$ {\bf on} $X$ is
\begin{eqnarray*}
{\rm Av}^w(F)(\gl)  = \lim_{R\to\infty}\frac{1}{\Vol(C_R)} \sum_{x\in C_R}\gl^w(\{x\}) T_xF(\gl) \\
{\rm Av}^w(F)(\gL) = \lim_{R\to\infty}\frac{1}{\Vol(C_R)} \sum_{x\in C_R}\gl^w(\{x\}) F(-x +\gL) \, ,
\end{eqnarray*}
if it exists. 

Prop.~\ref{averageThm} becomes: 
\begin{prop}\label{colouraverageThm}
The $w$-average value of  $F\in C(X)$ is defined at 
$\gL \in X $, $\mu$-almost surely and is almost surely equal to  $\palm^w(F)$. If $\mu$ is uniquely ergodic then the average value always exists and is equal to $\palm^w(F)$.
\end{prop}

We now come to the $w$-weighted autocorrelation. This is the measure on  $\RR^d$ defined by
\begin{align} \label{wACequation}
\gamma^{w}_{\gl}(f) & =  
\lim_{R\to\infty} \frac{1}{\Vol(C_R)}\gl^w|_{C_R} * \widetilde{\gl^w}|_{C_R}(f)\nonumber\\
 &= \lim_{R\to\infty} \frac{1}{\Vol(C_R)}\sum_{x \in C_R, y\in \Rd}\gl^w(\{x\}) 
 \overline{\gl^{w}(\{y\})}f(y-x)\\\nonumber
 &=  \lim_{R\to\infty} \frac{1}{\Vol(C_R)} \sum_{x \in C_R} 
 \gl^w(\{x\})N_f^{w}(T_{-x}\gl)\\  \nonumber
 &= {\rm Av}^w(N_f^{w}) =\palm^{w}(N_f^{w}) =: \palm^{w}_1(f) \nonumber
\end{align}
for all $f \in C_c(\RR^d)$ and for $\mu$-almost all $\gl \in X$.

We call $\palm^{w}_1$ the {\bf weighted first moment of the weighted Palm measure}.

\begin{theorem}\label{colourNautoCorr}  The weighted first moment $\palm^w_1$ of the weighted Palm measure is a positive definite measure. 
It is Fourier transformable and its Fourier transform $\widehat{\palm^w_1}$ is a
positive translation bounded measure on $\Rd$. Furthermore, $\mu$-almost surely, $\gl\in X$ admits a $w$-weighted autocorrelation $\gamma^w_{\gl}$ and it is equal to $\palm^w_1$.
If $X$ is uniquely ergodic then  $\palm^w_1 = \gamma^w_{\gl}$ for all $\gl \in X$.\qed
\end{theorem}  

\remark Regarding the statements about the transformability and translation boundedness
of the  Fourier transform, this is a consequence of the positive definiteness of
the Palm measure, see \cite{BF} Thm.~4.7, Prop.~4.9.

\smallskip

Prop.~\ref{link1} has the weighted form:
Let $g,h\in BM_c(\RR^d)$ and suppose that 
$g*\tilde h *\palm_1^w$ is a continuous function on $\RR^d$.
Then for all $t\in \RR^d$,
\begin{equation} \label{embeddingEq}
g*\tilde h*\palm_1^w (-t) = (T_t N^w_g, N^w_h) \, . 
\end{equation} 

Our interest now shifts to $L^2(\Rd, \widehat{\palm^w_1})$, its
inner product $\langle \cdot , \cdot \rangle^w$, and the unitary representation
$U^w$ of $\Rd$ on it which is given by the same formula as
(\ref{Urep}).

\subsection{The embedding theorem}
From equation (\ref{embeddingEq}) we obtain our embedding theorem, which is 
the full colour version of Prop.~\ref{link2}. 

\begin{theorem} \label{main}
For each system of weights $w=(w_1, \dots, w_m)$, the mapping
\begin{equation} \label{embeddingMap}
 \hat f \mapsto N^w_f \, ,
 \end{equation}
defined for all $f\in \mathbb{S}$, extends uniquely to 
an isometric embedding
\[\theta^w :L^2(\RR^d, \widehat{\palm_1^w}) \longrightarrow L^2(X, \mu) \]
which intertwines the represetations $U$ and $T$. 
\end{theorem}

We note here that the space on the left-hand side depends on $w$ while
the space on the right-hand side does not. The question of the image of
$\theta^w$ is then an interesting one. We come to this later.

We also note that the formula for $\theta^w(f)$ in (\ref{embeddingMap}), though true for $f \in \mathbb{S}$,
and no doubt many other functions too, is not true in general, and in particularly
not true for some functions that we will need to consider in the discussion of
spectral properties, e.g. see Cor.~\ref{spectrum}. 

Theorem \ref{main} gives an isometric embedding of $L^2(\Rd, \widehat{\palm^w_1})$ into $L^2(X,\mu)$ and along with it a correspondence of the spectral components of
$L^2(\Rd, \widehat{\palm^w_1})$ and its image in $L^2({X},\mu)$. Now the point is that
the spectral information of $L^2(\Rd, \widehat{\palm^w_1})$ can be read directly off that of
the measure $\widehat{\palm^w_1}$.
Specifically, let $\widehat{\palm^w_1} = (\widehat{\palm^w_1})_{pp} + (\widehat{\palm^w_1})_{sc}  + (\widehat{\palm^w_1})_{ac}$
be the decomposition of $\widehat{\palm^w_1}$ into its  pure point, singular continuous, and absolutely continuous parts. For $f \in L^2(\Rd, \widehat{\palm^w_1})$, the associated spectral measure $\sigma^w_f$ on $\Rd$ is given by
\[\langle f, U_t f \rangle^w = \int e^{2 \pi i x.t} d \sigma^w_f(x) \, .\]
However,
\[ \langle f, U_t f \rangle^w = \int e^{2 \pi i x.t} f(x)\overline{f(x)} d \widehat{\palm^w_1}(x) \, ,\]
so we have
\begin{equation} \label{sigmaf}
\sigma^w_f = |f|^2\widehat{\palm^w_1} = |f|^2(\widehat{\palm^w_1})_{pp} + |f|^2(\widehat{\palm^w_1})_{sc}  + |f|^2(\widehat{\palm^w_1})_{ac} \, ,
\end{equation}
which is the spectral decomposition of the  measure $\sigma_f$. With $\square$ standing
for pp, sc, or ac, we have
\begin{eqnarray*}L^2(\Rd, \widehat{\palm^w_1})_{\square} &:=& \{f \in L^2(\Rd, \widehat{\palm^w_1}) \;:\; \sigma^w_f \;\mbox{is of type}\; \square \} \\
&=& \{ f \in L^2(\Rd, \widehat{\palm^w_1}) \;:\; \supp (f) \subset \supp((\widehat{\palm^w_1})_{\square}) \} \, .
\end{eqnarray*}

This explains how information about the spectrum of the diffraction
can be inferred from the nature of the dynamical spectrum and
vice-versa. Since the mapping $\theta$ depends on $w$ and is not
always surjective, the correspondence between the two has to be
treated with care. Some examples of what can happen are given in
$\S$~\ref{Examples}.

Combining (\ref{sigmaf}) with Theorem~\ref{main}, we have S.~Dworkin's theorem:

\begin{coro}
Let $f\in C_c(\Rd)$. Then for $\mu$- almost all  $\gl \in X$,
$\widehat{ \gamma^w_{f*\gl}}$ is  the spectral measure $\sigma_{N^w_f}$ on
$L^2(X, \mu)$.
\end{coro}

 {\sc \bf Proof:} $\gamma^w_{f*\gl} = f*\tilde f *\gamma^w_{\gl}$,
 so $\widehat{\gamma^w_{f*\gl}} = |\widehat f|^2 \widehat{\gamma^w_{\gl}} =
 |\widehat f|^2 \widehat{\palm^w_1} = \sigma^w_{\widehat f}$ almost surely.
 Now,
 \[ \langle \widehat f, U_t \widehat f\rangle^w = (N^w_f, T_t N^w_f)_{L^2(X, \mu)} \]
so the spectral measure $\sigma^w_{\widehat f}$ computed for $L^2(\Rd, \widehat{\palm^w_1})$ is the same as the spectral measure $\sigma_{N^w_f}$ computed for $L^2(X, \mu)$. \qed

\begin{coro} \label{acscpp} 
For all $f, g \in L^2(\Rd, \widehat{\palm^w_1})$, the spectral measures
$({\langle U_t f,g\rangle^w})^{\vee}$ and $(T_t \theta^w(f), \theta^w(g))\iFT$ on $\Rd$ are equal,
and in particular of the
same spectral type:
absolutely continuous, singular continuous, pure point, or mix of these.
\end{coro}

\begin{coro}\label{spectrum}
For $k \in \Rd$, $\chi_k$ is in the point spectrum\footnote{One often simply says that
$k$ is in the point spectrum, with the understanding that it means $\chi_k$.} of $U_t$ if and only if
$\widehat{\palm^w_1}({k}) \ne 0$. The corresponding eigenfunction is $\chF_{\{-k\}}$. 
When this holds, $\chi_k$ is in the point spectrum of $T_t$,
the eigenfunction corresponding to it is $\theta^w(\chF_{\{-k\}})$, and 
$||\theta^w(\chF_{\{-k\}})|| = \widehat{\palm^w_1}(k)^{1/2}$.
\end{coro}

{\bf \sc Proof:} The first statement is clear from (\ref{Urep}) and our remarks above. For
the second, suppose that $k\in\Rd$ and $\widehat{\palm^w_1}({k}) \ne 0$.
Let $f\in L^1(\Rd, \widehat{\palm^w_1})$ be an eigenfunction for $k$. Then
 \[ \exp(2 \pi i k.t) f(x)= U_t f(x) = \exp(-2 \pi i t.x) f(x)\]
 for all $x \in \Rd$. For $x$ with
$f(x) \ne 0$, $\exp(2 \pi i k.t) = \exp(-2 \pi i x.t)$ for all $t \in \Rd$, so
$x= -k$. Thus $f = f(-k)\chF_{\{-k\}}$.
By Thm.~\ref{main}, $\theta^w(f)
\in L^2(X,\mu)$ with $T_t(\theta^w (f)) = \chi_{k}(t) \theta^w(f)$
for all $t\in \Rd$. 
 \qed
 
 \remark: One should note that the eigenvalues always occur in 
 pairs $\pm k$ since $\palm_1$ is positive-definite and
$\widehat{\palm^w_1}({-k})= \widehat{\palm^w_1}({k})$ . How does one work out $\theta^w(\chF_{\{-k\}})$? This is the
 content of the $L^2$-mean form Bombieri-Taylor conjecture that we shall establish in Sec.~\ref{SecBombieri-Taylor}.

\section{The algebra generated by the image of $\theta$} \label{algebra}

\subsection{The density of $\Theta^w(\mathbb{S})$}

\begin{theorem} \label{AlgebraGenerates}
Let $(X, \mu)$ be an $m$-coloured stationary uniformly discrete
ergodic point process and $w$ a system of weights. 
Suppose that the weights $w_i$, $i = 1,\dots,m$ are all different from
one another and also none of them is equal to $0$.  Then the algebra $\Theta^w$ 
generated by $\theta^w(\mathbb{S})$
and the identity function $1_X$ is dense in $L^2(X, \mu)$.
\end{theorem}
\begin{remark}
If $\widehat{\palm^w_1}({0}) \ne 0$
then $\theta^w(\mathbb{S})$ already contains $1_X$ by Cor.~\ref{spectrum}.
\end{remark}

The remainder of this subsection is devoted to the proof of this theorem. 

We begin with the construction of certain basic types of finite partitions of $X$.
Here we will find it easier to deal with coloured point
sets than with their corresponding measures.

Let $r>0$ be fixed so that $X\subset M_p(C^{(m)}_r,1)$. For each pair
of measurable sets $K,V \subset \Rd$, with $K$ bounded and
$V$ a neighbourhood of $0$, we define

\begin{equation} \label{localTop'}
U(K, V) := \{(\gL,\gL') \in \cDr \, : \,   K \,\cap  \,\gL \subset V + \gL' \quad \mbox{and} \quad  K \,\cap  \,\gL' \subset V + \gL \} \,,
\end{equation}
which is just a variation on (\ref{localTop}), and serves to define another
fundamental system of entourages for the same
uniformity, and then the same topology, on $X$ as we have been using all along. 
For any $\Phi \in \cDr^{(m)}$ we define
\[ U(K,V)[\Phi] := \{ \gL \in X : (\gL, \Phi) \in U(K,V) \}  \, .
 \]

We begin by choosing a finite grid in $\Rd$ and partitioning $X$ according
to the colour patterns it makes in this grid. Here are the details.
Let $K\subset \Rd$ be a half open cube of the form
 $[a_1,a_1 +R) \times \dots \times [a_d,a_d +R)$, $R>0$, and $V$ be an half-open
cube of
diameter less than $r$, centred on $0$, which is so sized that its translates
can tile $K$ without overlaps. The set of translation vectors used to make up
this tiling is denoted by $\Psi$, so in fact this set is the set of centres of
the tiles of the tiling. Each centre locates a tile and in each of these
tiles we can have at most one coloured point of $\gL$, that is, 
at most one pair $(x,i)$  with $x\in \Rd$ and $i\le m$. Let
\begin{alignat}{1}
\FrP :=\{\Phi= (\Phi_1, \dots, \Phi_m) \,:\, &(\Phi_0, \Phi_1, \dots, \Phi_m)\\
&\mbox{is an ordered partition of} \; \Psi \}  \, ;\nonumber
\end{alignat}
that is, we take all possible ordered partitions of $\Psi$ into $m+1$ pieces, which we interpret as all the various coloured patterns of cells of our tiling. 
$\Phi_i$ designates the cells containing the points of colour $i$ (second component $i$), $i =1,\dots m$ , and $\Phi_0$ designates all the cells which contain no points of the pattern. 

The inclusion relation $\subset$ on $\FrP$
by $\Phi= (\Phi_1, \dots, \Phi_m) \subset \Phi'= (\Phi_1', \dots, \Phi_m')$ if and only if $\Phi_i \subset \Phi_i'$, for all $1\le i \le m$, provides
a natural partial ordering on $\FrP$.
Using the  notation established in (\ref{basicOps}), for each  $\Phi \in \FrP$ define
 \[P[\Phi] := \{\gL \in X :  K\cap \gL  \subset V+ \Phi ,  K \cap \gL  \nsubseteq
 V+\Phi'
 \,\textrm{for any} \, \Phi' \varsubsetneqq \Phi \} \, . \]
 Because of the choice of $V$, an element of $X$ can have at most one point
in any one
 of the cubes making up the tiling of $K$. Each $P[\Phi]$ is the set of
 elements of $X$ which make the coloured pattern $\Phi$ inside the cube $K$.

 \begin{lemma}
 \[X = \bigcup_{\Phi \in \FrP} P[\Phi] \]
 is a partition of $X$. Furthermore, for all $\Phi \in \FrP$,
 \[U(K,V^\circ)[\Phi] \cap X \subset P[\Phi] \subset U(K,\overline
V)[\Phi] \cap X \, .\]
 \end{lemma}
 {\sc Proof:} By construction the $P[\Phi]$ form a partition of $X$. Let $\Phi \in
\FrP$ and let $\gL \in P[\Phi]$. Then
 $K \cap \gL  \subset V+ \Phi $. Also, for each
 $s$ lying in some component $\Phi_i$ of $\Phi$ there is $x\in K\cap \gL_i $ with $x \in V+ s$,
 whence $ s \in -V + x \subset \overline V +x$.
 This shows $K \cap \Phi \subset \overline V + \gL$, so
 $\gL \in U(K,\overline V)[\Phi]$.

 On the other hand, if $\gL \in U(K,V^\circ)[\Phi]$ then $K \cap \gL 
\subset  V^\circ + \Phi \subset
 V+ \Phi$, which is the first condition for $\gL \subset P[\Phi]$. Since also
$\Phi =  K\cap \Phi  \subset
 V^\circ + \gL $, for each $s$ in some component $\Phi_i$ there is $x \in \gL_i$ with
 $x = -v +s \in V^\circ +s  \subset V +s$.
 By the construction of the tiling of $K$, no other set $V^\circ + t$,
 $t \in \Psi$,  can contain $x$. Thus $\gL$
 meets every tile centred on a point of $\Phi$ and $\gL \in P[\Phi]$. \qed
 
 \smallskip
 
We know that $\theta^w(\mathbb{S})$ contains all the functions
$N^w_f$, $f\in \mathbb{S}$, in particular all the $N^w_f$, $f\in \mathbb{S}_c$, and so its $L^2$-closure contains $N^w_f$, $f \in C_c(\Rd)$ (use Prop.~\ref{convInNf}). Again, using 
Prop.~\ref{convInNf} we can conclude that $\overline{\theta^w(\mathbb{S})}$ contains
all the functions $N^w_A$, where $A$ is an bounded open or closed subset
of $\Rd$. We start with these functions and work to produce more complicated ones.

 \begin{lemma}\label{simple2}
  Let $s \in \Psi$ and let $i \le m$.
  Then the functions $N_{V^\circ +s}\circ {\rm res}^i(\gL)$ and
   $N_{\overline{V}+s}\circ {\rm res}^i(\gL)$ are in $\overline{\Theta^w}$.
  \end{lemma}

 {\sc Proof:} $N^w_{V^\circ+ s} \in \overline{\Theta^w}$. From (\ref{useOfres}) and  ${\textrm {diam}}(V) <r$,
\begin{eqnarray}
N^w_{V^\circ +s}(\gL)
&=& \sum_{i=1}^m w_i N_{V^\circ +s} \circ {\rm res}^i(\gL) \\
&=& \sum_{i=1}^m w_i N_{V^\circ +s}(\gL^i)  =  0\; \textrm{or} \; w_{j}
\end{eqnarray}
according as $(V^\circ +s) \cap \gL $ is empty or contains a (necessarily unique) point $x$ of some colour $j$. Write $\uF$ for
$N^w_{V^\circ +s}$ and $F$ for $N_{V^\circ +s}$. The first is a function on
$X$, the second a function on $r$-uniformly discrete subsets of $\Rd$ (see
(\ref{countingFunctions})).
 Then\footnote{Here the superscripts really mean powers!}
$\uF^j(\gL) = \sum_{i=1}^m w_i^j F(\gL_i)$ since always $F^j(\gL_i)= F(\gL_i)$
and $F(\gL_i)F(\gL_k) = 0$
whenever $i\ne k$. 

Let $W$ be the $m \times m$ matrix defined by $W_{jk} = w^j_k$, $1\le j,k \le m$. By the  hypotheses
on the weights it has an inverse $Y$. Then
\[ \sum_{j=1}^m Y_{ij} \uF^{j}(\gL) = \sum_{j=1}^m Y_{ij}\sum_{k=1}^m w_k^j F(\gL_k)
= F(\gL_i) \,. \]

This proves that the functions
$\gL \mapsto N_{V^\circ +s}(\gL_i)= N_{V^\circ +s}\circ {\rm res}_i(\gL)$
are all in $\Theta^w$.
The same argument applies in the case of $\overline V$. \qed

\begin{lemma}
For all $\Phi \in \FrP$, $\chF_{P[\Phi]} \in \overline{\Theta^w}$.
\end{lemma}

 {\sc Proof:}
Let $\Phi \in \FrP$ and assume $\Phi \ne \emptyset$.
  Let
  \begin{eqnarray*}
 f_1 :&=& \Pi_{i=1}^m \Pi_{s \in \Phi_i} N_{V^\circ +s}\circ {\rm res}_i \, \\
 f_2 :&=&\Pi_{i=1}^m \Pi_{s \in \Phi_i} N_{\overline{V} +s}\circ {\rm res}_i \, ,
  \end{eqnarray*}
  are all in $\Theta^w$. These functions take the value $1$ only on sets $\gL$ which
  hit all the cells $V^\circ +s$ (respectively $\overline{V} +s $)
  centred on the points and with the colours specified
  by $\Phi$. However, such $\gL$ may hit other cells also, hence
    \[f_1 \le \sum_{\Phi\subset\Phi' \in \FrP} \chF_{P[\Phi']} \le f_2  \, . \]
  However, for any fixed $i$, 
  \begin{eqnarray*}
  \int | N_{\overline{V} +s}(\gL_i) &-& N_{V^\circ +s}(\gL_i)|^2 \dmu(\gL) = 
  \int | N_{\overline{V} +s}(\gL_i) - N_{V^\circ +s}(\gL_i)| \dmu(\gL) \\
  &=& \int | N_{(\overline{V}\backslash V^\circ) +s}(\gL_i) | \dmu(\gL)  =
  I^i \Vol((\overline{V}\backslash V^\circ +s)) = 0
  \end{eqnarray*}
 showing that $N_{\overline{V} +s}$ and  $ N_{V^\circ +s}$ are equal as $L^2$
functions,
  whence also $f_1$ and $f_2$ are equal.
 This shows that 
 \[
 \sum_{\Phi'\supset\Phi}\chF_{P[\Phi']} \in \overline{\Theta^w} \, . \]

 In the case that $\Phi$ is empty,

\[\sum_{\Phi'\supset\Phi}\chF_{P[\Phi']} = \sum_{\Phi' \subset \Psi}
\chF_{P[\Phi']}  = \chF_{X} \in \overline{\Theta^w} \, . \]

 Now by  M\"obius inversion on the partially ordered on the subsets of $\FrP$ ,
 $\chF_{P[\Phi]} \in \overline{\Theta^w}$ for all $\Phi \in \FrP$.
  \qed

 \begin{lemma}\label{simple3}
 Let $F:X \longrightarrow \RR$ be a continuous function and let $\epsilon
>0$. Then there
 exist half-open cubes $K,V$ as above so that for the corresponding partition of
$X$,
 \[ || F - \sum_{\Phi \in \FrP} m_{\Phi} \chF_{P[\Phi]} ||_\infty \le
\epsilon \, ,\]
 where $m_{\Phi} := \inf \{ F(\gL) : \gL \subset P[\Phi]\}$.
 \end{lemma}
 {\sc Proof:} Since $X$ is compact, $F$ is uniformly continuous. Then
 given $\epsilon >0$ there
exist a compact set
 $K \subset \Rd$ and a neighbourhood $V'$ of $0 \in \Rd$ so that
 $|F(\gL) - F(\gL')| < \epsilon$ for all $(\gL,\gL') \in U(K,V') \cap (X
\times X)$. We can increase $K$ to
 some half-open cube of the type above without spoiling this and then choose some half-open cube
 $V$,
centred on $0$
 and of diameter
 less than $r$, which tiles $K$ and also satisfies $2\overline V \subset V'$. We let
 $\FrP$ be the corresponding set of partitions.

 Let $\Phi \in \FrP$. If $\gL,\gL'
\in
 U(K,\overline V)[\Phi] \cap X$ then $K \cap \gL  \subset \overline
V+\Phi $
 and $\Phi \subset \overline V+ \gL' $.  Thus for any $x\in K \cap \gL $, $x = v +s
= v'+ v+ x' $
 where  $s \in \Phi$, $x'\in \gL'$ (both with the same colour as $x$), and  $v,v'\in \overline V$, from which we conclude
 $K\cap \gL \subset 2\overline V + \gL'$. In the same way  $K \cap \gL' \subset
2\overline V + \gL$,
 so $(\gL, \gL') \in U(K,V')$ and $|F(\gL) - F(\gL')| < \epsilon$. In
particular this
 holds for all $\gL,\gL' \in P[\Phi]$, since it is contained in 
 $U(K,\overline V)[\Phi]$,
 and so $F$ varies by less than
 $\epsilon$ on $P[\Phi]$ . The result follows at once from this. \qed

\smallskip

The proof of Thm.~\ref{AlgebraGenerates} is an immediate consequence of this.
$\overline{\Theta^w(\mathbb{S})}$ contains the functions$\chF_{P[\Phi]}$ and so
also all their limit points, and hence all continuous functions on $X$.
Finally the continuous functions are dense in $L^2(X,\mu)$. \qed

 \smallskip
 
For the case $m=1$, recall that the $n$th moment of $\mu$ is the measure $\mu_n$ on $(\Rd)^n$ is defined by 
$\mu_n(f_1,\dots, f_n) = \mu(N_{f_1} \dots N_{f_n})$. Since Thm.~\ref{AlgebraGenerates} says that the linear
span of all the product functions $N_{f_1} \dots N_{f_n}$ is dense in $L^2(X,\mu)$,
we see that $\mu$ is entirely determined by its moment measures.

In the general case we may define the $nth$ weighted moments by:
\begin{equation}
\mu^w_n(f_1,\dots, f_n) = \mu(N^w_{f_1} \dots N^w_{f_n}) \,.
\end{equation}
For example, using \eqref{colourIntensity} and \eqref{weightedIntensity},
\begin{eqnarray} \label{firstWeightedMoment}
\mu^w_1(\chF_A) &=& \mu(N^w_A) = \int_X \sum w_i \gl^i(A) d \mu(\gl) \\
&=& \sum w_i I^{(i)} \Vol (A)  = I^w \Vol(A)\, ,\nonumber
\end{eqnarray} 
which shows that $\mu^w_1 = I^w\Vol$ (also see \eqref{intensityFormula}).

Then the same argument as in Thm.~\ref{AlgebraGenerates} leads to:

\begin{prop} \label{momentsToLaw} Let $(X, \mu)$ be an $m$-coloured stationary uniformly discrete
ergodic point process and $w$ a system of weights in which $w_i$, $i = 1,\dots,m$ are all different from one another and also none of them is equal to $0$. 
Then the measure $\mu$ is determined entirely by its set of $nth$ weighted
moments, $n = 1,2, \dots$ . 
\end{prop} \qed

We will relate this to higher correlations in the next section.

\begin{coro}
Let $(X,\mu)$ and $w$ be as in Prop.~{\rm \ref{momentsToLaw}}. Then the measure 
$\widehat{\palm^w_1}$ (which is the almost sure diffraction for
the members of $X$ when the weighting is $w$) is pure point if and only if the dynamical system $(X,\mu)$ is pure point, i.e., the
linear span of the eigenfunctions is dense in $L^2(X,\mu)$.
\end{coro}

\begin{remark}
This is the principal result of ~\cite{LMS1}. See also \cite{Gouere}.
\end{remark}

{\sc Proof:} The `if' direction is a consequence of Corollary~\ref{acscpp} of
Theorem~\ref{main}.

The idea behind the `only if' direction is simple enough. The assumption
is that the linear space of the eigenfuctions $L^2(\Rd,\widehat{\palm^w_1})$
is dense in $L^2(\Rd,\widehat{\palm^w_1})$, and eigenfuctions of 
$L^2(\Rd,\widehat{\palm^w_1})$ map to 
eigenfuctions of $L^2(X,\mu)$ under $\theta^w$. However, products of eigenfunctions of $(X,\mu)$ are again eigenfunctions. We know that the algebra generated by
the image of $\mathbb{S}(\Rd)$ in $L^2(X,\mu)$ is dense. So the linear space
that we get by taking the algebra generated by the eigenfuctions ought also to be both dense and linearly generated by eigenfunctions. The trouble is that the eigenfuctions of $L^2(\Rd,\widehat{\palm^w_1})$ are not in $\mathbb{S}(\Rd)$ and the 
space $L^2(X,\mu)$ is not closed under multiplication, so we need to be careful. 

The set $BL^2(X,\mu)$ of measurable square integrable functions on $X$ that 
are bounded on a subset of full measure form an algebra (i.e. the product of such functions are also bounded), and $\theta^w(\mathbb{S})$ is contained in it. In fact
any bounded function of $L^2(\Rd,\widehat{\palm^w_1})$ is mapped by 
$\theta^w$ into $BL^2(X,\mu)$, as we can see from Theorem~\ref{main}
and equation (\ref{weightedN}) and taking approximations by elements of
$\mathbb{S}$. For $F \subset BL^2(X,\mu)$, let $L(F)$ denote its linear span
and $\langle F \rangle_{\rm alg}$ the subalgebra of $BL^2(X,\mu)$
generated by $F$.

By Corollary 3, $\chi_{k}$ is in the point spectrum of
$U_t$ if and only if $\widehat{\palm^w_1}(k)\neq 0$, and the eigenfunction
corresponding to $\chi_{k}$ is $\theta^w(\chF_{\{-k\}})$.
Denote by $E$ the set of
$\{\chF_{\{-k\}}:\widehat{\palm^w_1}(k)\neq 0\}$ and by 
$L(E)$ its linear span. By Theorem~\ref{main}, $\theta^w(E)$ is a set of
eigenfunctions of $T_t$, and by what we just saw
$\theta^w(L(E))\subset BL^2(X,\mu)$. By assumption, $L(E)$ is
dense in $L^2(\Rd,\widehat{\palm^w_1})$.

Then
\[\overline{L(\theta^w(E))} \supset \theta^w(\overline{L(E)}) 
\supset \theta^w(\mathbb{S}) \, \]
and 
\[BL^2(X,\mu) \supset \langle \theta^w(L(E))\rangle_{\rm alg} \, .\]
Thus,
\[ \overline{\langle \theta^w(E) \rangle_{\rm alg} }
= \overline{\langle \theta^w(\overline{L(E)}) \rangle_{\rm alg} } 
\supset \overline{\langle\theta^w(\mathbb{S})\rangle_{\rm alg}} = L^2(X,\mu) \, ,\]
which shows that the denseness of the linear span of the eigenfunctions of $L^2(X,\mu)$.
\qed

\section{Higher correlations and higher moments}

Let $\xi: (\Omega, \CalA, P) \longrightarrow (X,\CalX)$ be a uniformly discrete
stationary multi-variate point process with accompanying dynamical 
system $(X,\RR^d, \mu)$, and let $w$ be a system of weights.

The $n+1$-point {\bf correlation} ($n = 1,2, \dots$) of $\gl \in X$ is the measure
on $(\Rd)^n$ defined by 
\begin{eqnarray*}
\gamma_\gl^{(n+1)}(f)&=&\lim_{R\to\infty} \frac{1}{\Vol (C_R)}\sum_{y_1,\dots, y_n,x
\in C_R} 
\gl^w(\{x\})\Pi_{i=1}^n\gl^w(\{y_i\}) T_xf(y_1, \dots,y_n)\\
&=&\lim_{R\to\infty} \frac{1}{\Vol (C_R)}\sum_{\stackrel{x\in C_R}{ y_1,\dots y_n \in \Rd}} 
\gl^w(\{x\})\Pi_{i=1}^n\gl^w(\{y_i\}) T_xf( y_1, \dots, y_n) \, ,
\end{eqnarray*}
for all $f\in C_c((\Rd)^n)$. 
In particular for $f =(f_1, \dots, f_n) \in (C_c(\Rd))^n$, where each 
$f_i \in C_c(\Rd)$, 
\begin{eqnarray*}
\gamma_\gl^{(n+1)}(f)
&=&\lim_{R\to\infty} \frac{1}{\Vol (C_R)}
\sum_{\stackrel{x\in C_R}{ y_1,\dots y_n \in \Rd}}
\gl^w(\{x\})\Pi_{i=1}^n\gl^w(\{y_i\}) T_xf_1(y_1) \dots  T_xf_n(y_n)\\
&=&\lim_{R\to\infty} \frac{1}{\Vol (C_R)}\sum_{x\in C_R} 
\gl^w(\{x\})N^w_{f_1}(T_{-x} \gl) \dots  N^w_{f_n}(T_{-x}\gl )\\
&=& \Av^w( N^w_{f_1}\dots  N^w_{f_n} )(\gl) \, .
\end{eqnarray*}

We know that $\mu$-almost surely this exists, it is independent of $\gl$, and
\[  \Av^w( N^w_{f_1}\dots  N^w_{f_n} )(\gl) = \palm^w(N^w_{f_1}\dots  N^w_{f_n})
=: \palm^w_n( (f_1, \dots , f_n))\,.\]
The measure defined on  the right-hand side of this equation is the $n$th 
{\bf weighted moment of
the weighted Palm measure} at $(f_1, \dots f_n)$, so we arrive at the useful fact which 
generalizes what we already know for the $2$-point correlation:
\begin{prop} \label{corr=palmMoment}
The $n+1$-point correlation measure exists almost everywhere on $X$ and
is given by $\palm^w_n$. 
\end{prop}

Of course, in the one colour case where there are no weights (or if the weighting is trivial: $w = (1,\dots, 1)$), then these are ordinary moments.

\begin{lemma}\label{intensityTrick} Assume that the weights are all different and none of them is zero. Then the weighted intensity (and hence the 
first moment of $\mu$) is determined by the $\palm^w_n$, $n =1, 2, \dots$  \, .
\end{lemma}

\noindent
{\sc \bf Proof}: By \eqref{intensityFormula} we need to know $\palm^w(X)$. Now
$\palm^w$ is supported on the set $X_0$ of elements $\gl \in X$ which
have an atom at $\{0\}$ and we have
\begin{eqnarray*}
\palm^w(X)&=& \int_{X_0} N_{\chF_{\{0\}}}(\gl^{\downarrow})d\palm^w(\gl)
= \int_{X_0}\sum_{i=1}^m N_{\chF_{\{0\}}}(\gl^i)d\palm^w(\gl)\\
&=&\sum_{i=1}^m \palm^w(N_{\chF_{\{0\}}}\circ \mbox{res}^i).
\end{eqnarray*}
From $\palm^w_1$ we have
\[
\palm^w_1(\{0\})= \int_{X_0} N^w_{\chF_{\{0\}}}(\gl)d\palm^w(\gl)
=\sum_{i=1}^m w_i\palm^w(N_{\chF_{\{0\}}}\circ \mbox{res}^i) \, .
\]
Similarly,  
\[
\palm^w_2(\{(0,0)\})=\int_{X_0} (N^w_{\chF_{\{0\}}})^2(\gl)d\dot{\mu}^w(\gl)\\
=\sum_{i=1}^m w_i^2\dot{\mu}^w(N_{\chF_{\{0\}}}\circ \mbox{res}^i)\, .
\]
Continue this until we get to
\[
\palm^w_m(\{(0,\dots,0)\})=\int_{X_0} (N^w_{\chF_{\{0\}}})^m(\gl)d\palm^w(\gl)\\
=\sum_{i=1}^m w_i^m\dot{\mu}^w(N_{\chF_{\{0\}}}\circ \mbox{res}^i) \, .
\]
Using the same argument in Lemma~\ref{simple2}, we can solve this system
of equations for
$\palm^w(N_{\chF_{\{0\}}}\circ\mbox{res}^i)$ for $i=1,\dots,m$
and hence determine the weighted intensity $\palm^w(X)$. \qed

\begin{theorem}\label{correlationsDetermineMu}
Let $(X, \mu)$ be an $m$-coloured stationary uniformly discrete
ergodic point process and $w$ a system of weights in which $w_i$, 
$i = 1,\dots,m$, are all different from one another and none of them is equal to $0$. Then the measure $\mu$ is completely determined by 
the weighted $n+1$-point correlations of 
$\mu$-almost surely any $\gl \in X$, $n = 1,2, \dots$. 
\end{theorem}

The key to this is the known fact (in the non-weighted case) that the 
$n$th moment of the Palm measure, $n = 1,2, \dots$, is the same as the reduced 
$(n+1)$st moment of the measure $\mu$ itself. Thus knowledge of the 
correlations gives us the moments $\palm_n$  of the Palm
measure, which in turn is the same as knowledge of the reduced moments
of $\mu$. These in turn determine the moments $\mu_{n+1}$, $n=1,2, \dots$
of $\mu$. As for $\mu_1$,
we already know that it is just the intensity of the point process times Lebesgue measure, and from Lemma~\ref{intensityTrick}, this is derivable from the 
moments.  

First of all we give
a short derivation of these facts in the unweighted $m=1$ case, and then
show how to augment these to the weighted case.

Let $g,h_1, \dots h_n \in C_c(\Rd)$ be chosen freely. Let $G:\Rd \times X\longrightarrow \CC$ be defined by
\[G(x,\gl) = g(x) N_{T_xh_1}(\gl) \dots N_{T_xh_n}(\gl) \,. \]
We use the Campbell formula
\[\int_X\sum_{x\in\Rd} \gl(\{x\}) G(x,\gl) d\mu(\gl) 
=\int_{\Rd} \int_X G(x, T_x\gl) d\palm(\gl)  \,.\]
The left-hand side reads\footnote{
$g(T_xh_1) ... (T_x h_n)$ stands for the function whose value on
$(x,(y_1, \dots y_n)) \in \RR^d \times (\RR^d)^n$ is $g(x)(T_xh_1)(y_1) ... (T_x h_n)(y_n)$. }
\begin{alignat}{1} \label{reducing}
\int_X\sum_{x\in\Rd} &\gl(\{x\})g(x) N_{T_xh_1}(\gl) \dots N_{T_xh_n}(\gl)d\mu(\gl)\\\nonumber
&= \int_X \gl(g) \gl(T_xh_1) \dots \gl(T_xh_n) d\mu(\gl) \\ \nonumber
&= \mu_{n+1}(g (T_xh_1) \dots (T_xh_n)) 
= \int_{\Rd}g(x) dx \, \mu^{\rm red}_{n+1}(h_1\dots h_n) \, , \nonumber
\end{alignat}
while the right-hand side reads
\[ \int_{\Rd}\int_X g(x) N_{T_xh_1}(T_x\gl) \dots N_{T_xh_n}(T_x\gl)
d\palm(\gl) dx = \int_{\Rd}g(x) dx \, \palm_n(h_1\dots h_n) \, ,\]
since $N_{T_xh}(T_x\gl) = N_h(\gl)$.
For the reduced moments see \cite{DV-J}, Sec. 10.4, especially
Lemma 10.4.III and Prop. 10.4.V. The point is that
$\mu_{n+1}$ is invariant under simultaneous translation of its $n+1$ variables. This
invariance can be factored out leading to the rewriting of $\mu_{n+1}$ as a product of
Lebesgue measure and another measure, which is, {\em by definition}, the reduced measure.
This rewriting is exactly the last part of equation (\ref{reducing}).
Thus, $\mu^{\rm red}_{n+1} = \palm_n$ and,
using Prop.~\ref{momentsToLaw} and Prop.~\ref{corr=palmMoment}, 
Theorem~\ref{correlationsDetermineMu} is proved in the $1$-coloured case.

To obtain the weighted version, we use now the functions
\[G^w(x,\gl) = g(x) N^w_{T_xh_1}(\gl) \dots N^w_{T_xh_n}(\gl)  \]
and the weighted form of the Campbell formula. Then the same argument leads to
$(\mu^w_{n+1})^{\rm red} = \palm_n^w$, $n =1, 2, \dots $. Meanwhile, the first moment
is determined by the weighted intensity, given in Lemma~\ref{intensityTrick}.
The proof of Theorem~\ref{correlationsDetermineMu} now follows as in the unweighted
case. \qed

\section{Examples} \label{Examples}

In this section we offer examples that show a variety of ways in which the image of the diffraction appears in the dynamics and in particular how the weighting system influences the outcome. We begin with a general construction of $m$-coloured uniformly discrete ergodic point processes from symbolic shift systems, which allows one to lift results from the theory of the discrete dynamics of sequences to our situation of continuous dynamics.  

The first two of the examples come from well-known results about the 
Thue-Morse and Rudin-Shapiro sequences. Both sequences lead to dynamical systems 
$(X,\RR,\mu)$ of point sets on the real line, which are uniquely ergodic and minimal but for which the  mapping $\theta: L^2(\RR,\widehat{I_q}) \rightarrow L^2(X,\mu)$ is not surjective. In both cases the diffraction and dynamical spectra are mixed (pure point + singular in the one case, pure point + absolutely continuous in the other). However the mapping $\theta$ does not map the pure point diffraction surjectively to the pure point dynamical spectrum -- in fact, it can miss entire spectral components -- and this shows that $\theta$ itself is not in general surjective. This fact, that 
the diffraction does not convey full information on the dynamics has been pointed out much earlier by van Enter and Mi\c{e}kisz \cite{vEM}.

We then look at extinctions in model sets and observe that even in these most well-behaved sets, the diffraction and dynamical
spectra (both of which are pure point) need not match exactly. Finally we give an example to show the necessity of  non-zero weights in Theorem \ref{AlgebraGenerates}.

We begin with a short general review of the discrete dynamics of the sequences and look at what happens when we move to the continuous setting by using suspensions. We have done this in slightly more generality than we need for the examples, but with a view to further applications \cite{Deng}. 

\subsection{The continuous dynamics of sequences on the real line}

A good source for examples is to start with symbolic shifts. We start with a finite alphabet ${\bf m} =\{1,\dots m\}$ and then define
${\bf m}^\ZZ$ to be the set of all bi-infinite sequences 
$\zeta = \{z_i\}_{-\infty}^{\infty}$,
which we supply with the product topology. 
Along with the usual shift action $(T(\zeta))_i = \zeta_{i+1}$ for all $i$, ${\bf m}^\ZZ$ becomes a dynamical system over the group $\ZZ$. We are interested in compact $\ZZ$-invariant subspaces $X_\ZZ$ of $({\bf m}^\ZZ, \ZZ)$. We will assume that $(X_\ZZ,\ZZ)$ is equipped with an invariant and ergodic probability measure $\mu_\ZZ$. Such measures always exist. 
We define for all 
$\underline{i} = (i_0, \dots, i_k) \subset {\bf{m}}^{k+1}$, $k = 0,1 \dots $, and $p\in \ZZ$,
\[X_\ZZ[\underline{i};p] = \{\zeta \in X_\ZZ \,:\, z_{j +p} = i_j, j = 0,\dots k\} \,. \] 
These cylinder sets form a set of entourages for the standard uniform topology on
$X_\ZZ$ which defines the product topology. When $p =0$, we usually leave it out and also leave off the parentheses; so, for example,  $X_\ZZ[ij]$ means $X_\ZZ[(ij);0]$. 

We need to move from the discrete dynamics (action by $\ZZ$) of $(X_\ZZ,\ZZ)$ to continuous dynamics with an $\RR$-action. There is a standard way of doing this by creating the suspension flow of $(X_\ZZ,\ZZ)$, and this new dynamical system has a natural invariant and ergodic measure and so satisfies our conditions PPI, PPII, PPIII. 
Basically each bi-infinite sequence 
$\zeta$ of $(X_\ZZ,\ZZ)$ is converted into a bi-infinite sequence of coloured points on the real line with $z_0$ being located at $0$. The most obvious thing is to space out the other points of the sequence on the integers, so that $z_n$ ends up at position $n$. The result can be viewed as a tiling of the line with coloured tiles of length $1$, the colour of a tile being the colour of the left end point that defines it. However, there are good reasons to allow different colours to have different tile lengths. \footnote{Readers interested only in the examples germane to this paper
may ignore the introduction of different tile lengths that we introduce here.}

For this purpose we take any set $\CalL= \{L_1, \dots, L_m\}$ of positive numbers as the tile
lengths, with an overall scaling so that 

\[ \sum_{j=1}^m L_j \mu_\ZZ(X_\ZZ[j]) = 1 \, .\]
Let $r = \min\{L_1, \dots,L_m\}/2$.

Given $\zeta = \{z_n\}_{-\infty}^\infty \in X_\ZZ$, define the sequence $S=S(\zeta) = \{S_n\}_{-\infty}^\infty$
by $S_0 = 0$, $S_n = \sum_{j=0}^{n-1}  L_{z_j}$, if $n >0$, $S_n = -\sum_{j=n}^{-1}  L_{z_j}$ if
$n <0$.

Define 
\begin{eqnarray}\label{piL}
\pi^{\CalL}: \RR \times X_\ZZ &\longrightarrow& \cDr^{(m)}(\RR)\\
(t,\zeta) &\mapsto & \{(t+ S_n, z_n)\}_{-\infty}^\infty \, ,\nonumber
\end{eqnarray}
which ``locates'' the symbols of $\zeta$ along the line (including colour information)
so that the $n$th symbol occurs at $t+ S_n$. This simultaneously provides us
with a tiling of the line by line segments of lengths $\{L_{z_n}\}$. We let $X^\CalL_\RR:= \pi^\CalL(X_\ZZ) \subset
\cDr^{(m)}(\RR)$. Both $\RR \times X_\ZZ$ and $\cDr^{(m)}(\RR)$ have natural 
$\RR$-actions on them, and the mapping $\pi^{\CalL}$ is $\RR$-invariant. It is easy to
see that $\pi^{\CalL}$ is continuous.

Let $R$ be the equivalence relation on
$\RR \times X_\ZZ$ defined by transitive, symmetric, and reflexive extension of
$(t,\zeta) \equiv_{R} (t + L_{z_0}, T\zeta)$. Evidently pairs are $R$-equivalent if and
only if they have the same image under $\pi^\CalL$. In fact, 
$(t,\zeta)$ is $R$-equivalent to a  unique element of 
\[ F^\CalL:= \bigcup_{i=1}^m (-L_i , 0] \times X_\ZZ[i] \]
and the mapping $\pi^\CalL$ is injective on this set. Since $\overline{F^\CalL} =
\bigcup [-L_i , 0] \times X_\ZZ[i]$ is compact and $\pi^\CalL$ maps this set onto $X^\CalL_\RR$,
we see that $X^\CalL_\RR$ is compact and hence $(X^\CalL_\RR, \RR)$ is a topological dynamical system.  

\subsection{Measures on the suspension} We define a positive measure $\mu^\CalL$ on $X^\CalL_\RR$ by 
\[\mu^\CalL(B) := (\Vol \otimes \mu_\ZZ )((\pi^\CalL)^{-1}(B) \cap F^\CalL) \]
for all Borel subsets $B$ of $X^\CalL_\RR$. We observe that $\mu^\CalL$ is a probability
measure since $\mu^\CalL(X^\CalL_\RR) = (\Vol \otimes \mu_\ZZ )(F^\CalL) = \sum L_i \mu_\ZZ(X[i]) = 1$.

This is an $\RR$-invariant measure on $X^\CalL_\RR$. It suffices
to show the shift invariance for sets of the form $J \times C $ where $J$ is an interval in $(-L_i,0]$ and $C$ is
a measurable subset of some $X_\ZZ[i]$, since these sets generate the $\sigma$-algebra of all Borel subsets of $F^\CalL$. We show that shifting of $J$ by $s<0 $ leaves the measure invariant. It is sufficient to do this for 
$|s| < \min\{L_1, \dots, L_m\}$, since we can repeat the process if necessary to account for  larger $s$. If $s+J \subset (-L_i,0]$, then the invariance of $\Vol$ gives what we need immediately.
If $s+J \nsubseteq (-L_i,0]$ then we may break $J$ into two parts; the part which remains in the interval
and the part which moves out of it to the left. We can restrict our attention to the part that moves out and then assume that $(s+J) \cap (-L_i,0] = \emptyset$. Then we
bring $(s+J )\times C$ back into $F^\CalL$ by
writing $C = \bigcup_{j=1}^m C \cap X_\ZZ[ij]$ so that
\[(s+J) \times C \equiv_R \bigcup_{j=1}^m (L_i +s +J) \times T(C\cap X_\ZZ[ij]) \,. \]
The measure of this is $\sum_{j=1}^m \Vol(J)\mu_\ZZ( T(C\cap X_\ZZ[ij])) = \Vol(J) \mu_\ZZ(C)= (\Vol \otimes \mu_\ZZ)(J\times C)$,
which is what we wished to show. 

 If the original measure $\mu_\ZZ$ on $X_\ZZ$ is ergodic, then so is the measure $\mu^\CalL$.
One way to see this is to start with the case when $\CalL = \{1,\dots, 1\}$. In this case we shall denote the objects that we have constructed above with a superscript ${\bf 1}$ rather than $\CalL$. It is easy to see that 
$\mu^{\bf 1}$ is an ergodic measure on $X^{\bf 1}_\RR$ since the latter can be thought of as $X_\ZZ \times U(1)$, where $U(1)$ is the unit circle in $\CC$, with the action of $\RR$ being such that going clockwise around the circle once returns one to the same sequence in $X_\ZZ$ except shifted once. 

We can define a flow equivalence $\phi: X^{\bf 1}_\RR \longrightarrow X^\CalL_\RR$ 
in the following way.
For each $\zeta \in X_\ZZ$ define $f_\zeta :\RR \longrightarrow \RR$ by
\[
f_\zeta(t) = \left\{ \begin{array}{ll} 
 |S_{-(k-1)}| + L_{z_{-k}}(t -|S_{-(k-1)}|) &
\mbox{if $t\ge 0, k-1\le t < k$} \\
-S_{k-1} + L_{z_k}(t -S_{k-1}) &
\mbox{if $t\le 0, k-1\le |t| < k$} \,.
\end{array} 
\right.  \]
This is a strictly monotonic piece-wise linear continuous function which fixes $0$.  
Its intent is clear: if $(t,\zeta)$ is understood to represent the sequence $\zeta$ placed down
in equal step lengths of one unit starting with $z_0$ at $t$, then $(f_\zeta(t), \zeta)$ represents
the same sequence, now scaled to the new colour lengths $L_{z_j}$ where $0$ is the fixed point.

Thus define a mapping $\RR \times X_\ZZ \longrightarrow \RR \times X_\ZZ$ by
$(t,\zeta) \mapsto (f_\zeta(t), \zeta)$. This mapping factors through the equivalence
relations that define $X^{\bf 1}_\RR$ and $X^\CalL_\RR$ to give the mapping $\phi$ which is the flow equivalence that we have in mind. For $I\times X_\ZZ[\underline{u}]$, where $I \subset (-L_{u_0}, 0]$, 
\[ \phi^{-1}((I\times X_\ZZ[\underline{u}])^\sim) = \left(\frac{I}{L_{u_0}} \times X_\ZZ[\underline{u}]\right)^\sim \,, \]
where the equivalence relations are taken for $\CalL$ and for ${\bf 1}$ respectively. Furthermore, 
$\mu^\CalL((I\times X_\ZZ[\underline{u}])^\sim) = \Vol(I) \mu_\ZZ(X_\ZZ[\underline{u}])$ and
$\mu^{\bf 1}((I/L_{u_0}\times X_\ZZ[\underline{u}])^\sim) = \Vol(I/L_{u_0} ) \mu_\ZZ(X_\ZZ[\underline{u}])$.

Now, if $B$ is an $\RR$-invariant subset of $X^\CalL_\RR$ then $\phi^{-1}(B)$ is an $\RR$-invariant
subset of $X^{\bf 1}_\RR$, and so, assuming that $\mu_\ZZ$ is ergodic, $\phi^{-1}(B)$ has measure
$1$ or $0$. If the former, then for all $i\le m$, $\phi^{-1}(B) \cap ((-1, 0]\times X_\ZZ[i])$ has 
$\mu^{\bf 1}$-measure $\mu_\ZZ(X_\ZZ[i])$ from which $B \cap ((-L_i, 0]\times X_\ZZ[i])$ has measure
$L_i \mu_\ZZ(X_\ZZ[i])$, which shows that $B$ is of full  measure in $X^\CalL_\RR$. A similar argument
works for the measure $0$ case. This shows that $\mu^\CalL$ is ergodic.

\subsection{Spectral features of the suspension} At this point we have
arrived at the setting of this paper: $(X^\CalL_\RR,\RR, \mu^\CalL)$ is a dynamical system satisfying PPI, PPII, PPIII. Henceforth we shall assume that the set of lengths
$\CalL = \{ L_1, \dots, L_m\}$ is fixed, and drop them from the notation. We may weight the system by choosing any real vector $w = (w_1, \dots w_m)$
of weights and assigning weight $w_i$ to the colour $a_i$. According to
Prop.~\ref{colourNautoCorr} , the weighted first moment $\palm^w_1$ of the weighted Palm measure is almost everywhere the weighted autocorrelation of the point sets of $X_\RR$, and this is everywhere true if the system is uniquely ergodic. We will use the symbol $w$ to also denote the mapping ${\bf m} \longrightarrow \{w_1, \dots w_m\}$,
$w(i) = w_i$.

We now come to the autocorrelation. For the purposes of the examples, it is convenient
to have all tile lengths equal to $1$:  $L_j =1$ for all $j$, and we shall assume this for the
remainder of this section.

Now let $\zeta = \{z_i\}_{-\infty}^{\infty} \in X_\ZZ$. Its autocorrelation, assuming that it exists, is
\[\gamma_\zeta^{w,\ZZ} = \sum_{k\in\ZZ} \eta^w(k)\delta_k^\ZZ\; , \]
defined on $\ZZ$,  where
\[\eta^w(k) := \lim_{N\to\infty} \frac{1}{2N+1}\sum_{i=-N}^N w(z_i)w(z_{i+k}) \, .\]
Its autocorrelation $\gamma_\zeta^{w,\RR}$ when thought of as an element of  $X_\RR$ is defined
on $\RR$ and is given by
\[\gamma_\zeta^{w,\RR} = \sum_{k\in\ZZ} \eta^w(k)\delta_k^\RR\; , \]
with the same $\eta^w(k)$.

The difference is in the delta measures, which are defined
on $\ZZ$ and $\RR$ respectively.  Thus 
$\widehat{\gamma_\zeta^{w,\ZZ}}$ is a
 measure on $\TT := \RR/\ZZ$ while $\widehat{\gamma_\zeta^{w,\RR}}$ is a measure
 on $\RR$. The relationship between these two measures is simple: for $x\in \RR$ and
 $\dot x := x \mod \ZZ$,
 \[  \widehat{\delta_k^\ZZ}(\dot x) = e^{-2 \pi i k.x} ,\quad
\widehat{\delta_k^\RR}(x) = e^{-2 \pi i k.x}   \, .\]
Thus, for all $k\in \ZZ$, $\widehat{\delta_k^{w,\RR}}$ is just the natural periodic extension of
$\widehat{\delta_k^{w,\ZZ}}$  and $\widehat{\gamma_\zeta^{w,\RR}}$ is the periodization of
$\widehat{\gamma_\zeta^{w,\ZZ}}$:
\[  \widehat{\gamma_\zeta^{w,\RR}}(x) = \widehat{\gamma_\zeta^{w,\ZZ}}(\dot x) \, .\]
The latter, hence also the former, exists almost surely.

The pure point, singular continuous, and absolutely continuous parts are also periodized
in this process and retain the same types. Thus if the pure point part of
$\widehat{\gamma_\zeta^{w,\ZZ}}$ is $\sum_{\dot k\in S} a_{\dot k}\delta_{\dot k}^\ZZ$ then
the pure point part of $\widehat{\gamma_\zeta^{w,\RR}}$ is $\sum_{k\in \RR,\, \dot k\in S} a_{\dot k}\delta_{k}^\RR$, where $\dot k = k \mod \ZZ$.

When it comes to $L^2(X_\ZZ,\mu_\ZZ)$ and $L^2(X_\RR, \mu)$ we make the following observation. If $f_{\dot k}$ is an eigenfunction for the action of $T$ on $L^2(X_\ZZ,\mu_\ZZ)$ corresponding to the eigenvalue $\dot k$ -- that is,
$T^nf_{\dot k} = \exp(2\pi i k.n) f_{\dot k}$  for some (or any) $k\in \ZZ$
with $\dot k = k \mod \ZZ$,
we can define a function $f_k$ on $X_\RR$ by
\[ f_k(t + \zeta) = \exp(-2 \pi i k.t) f_{\dot k}(\zeta)  \,. \]
It is easy to see that this is well-defined and is an eigenfunction for the
$\RR$-action on $X_\RR$ with eigenvalue $-k$ on $\RR$.
(The change in sign results from the fact that $T$ means shift left by $1$, whereas
$T_t$ means shift right by $t$.) This way we see that we have eigenfunctions for
$X_\RR$ which are all the possible continuous lifts of the eigenfunctions on $\RR/\ZZ$
to eigenfunctions on $\RR$.

Unfortunately there does not seem to be any simple connection between the 
other spectral components of $(X_\ZZ,\mu_\ZZ)$ and $(X_\RR,\mu)$. Thus,
for these components, we will
be reduced to the consequences that come by the embedding
of the diffraction into the dynamics.  

\subsection{The hull of a seqence} We start with an infinite sequence $\xi = (x_1, x_2, \dots)$
of elements of our finite alphabet ${\bf m}$ and define
$X_\ZZ(\xi)$ to be the set of all bi-infinite sequences $\zeta = \{z_i\}_{-\infty}^{\infty}
\in {\bf m}^\ZZ$ with the property that every finite subsequence $\{z_n, z_{n+1}, \dots, z_{n+k} \}$ ({\em word}) of $\zeta$ is also a word $\{x_p, x_{p+1}, \dots, x_{p+k}\}$
of $\xi$. 
Then set $X_\ZZ(\xi)$ is a closed, hence compact subset of ${\bf m}^\ZZ$, and $(X_\ZZ(\xi) ,\ZZ)$ is a dynamical system, called the {\em dynamical hull} of $\xi$. $(X_\ZZ(\xi),\ZZ)$ is mininal (every orbit is dense) if and only if $\xi$ is repetitive (every word reoccurs with bounded gaps).

Given a word $s = \{x_p, x_{p+1}, \dots, x_{p+k}\}$ of $\zeta$, we can ask about the frequency of its appearance (up to translation) in $\zeta$. Let
$L(s, [M,N])$ be the number of occurrences of $s$ in the
interval $[M,N]$.
The frequency of $s$ (relative to $t\in \ZZ$)  is  $\lim_{N\to\infty}L(s,t+[-N,N])/2N$,
if it exists.  It is known that the system $X_\ZZ(\xi)$ is 
both minimal and uniquely ergodic (that is, {\em strictly ergodic}) if and only if for for every $\zeta \in X_\ZZ(\xi)$ and every word $s$ of $\zeta$ the frequency of $s$  exists, the limit is approached uniformly for all in $t\in \ZZ$, and the frequency is positive. 
All of this is standard from the theory of sequences and symbolic dynamics
\cite{Quef}, Cor. IV.12 . 

We can transform $X_\ZZ(\xi)$ into a flow over $\RR$ by the technique discussed in the previous subsection and thus obtain $X_\RR(\xi)$, which will be minimal (respectively ergodic, uniquely ergodic) according as $X_\ZZ(\xi)$ is. 

In the next two subsections we consider situations which are derived from two famous sequences, the Thue-Morse and Rudin-Shapiro sequences.

\subsection{Thue-Morse}

The Thue-Morse sequence can be defined by iteration of the two
letter substitution (we use $a,b$ instead of $1,2$) 
\[a \rightarrow ab; \quad b\rightarrow ba: \quad
\xi = abbabaabbaababba \dots  \,. \]
based on the alphabet $A=\{a,b\}$ (we use $\{a,b\}$ instead of $\{1,2$\}).

Since the substitution is primitive, it is known that the
corresponding dynamical system $X_{\ZZ}=X_{\ZZ}(\xi)$, and hence
also $X_{\RR}=X_{\RR}(\xi)$, is minimal and uniquely ergodic.

For an arbitrary weighting system $w = (w_a,w_b)$ we have the diffraction
$w_a^2 \gamma_{aa} + w_aw_b \gamma_{ab} + w_bw_a \gamma_{ba} + w_b^2 \gamma_{bb}$ where $\gamma_{ij}$ is the correlation between points of types
$i,j \in A$. The natural symmetry $ a \leftrightarrow b$ of $X_{\ZZ}$ gives
$\gamma_{aa} = \gamma_{bb}, \gamma_{ab}= \gamma_{ba}$. 

Kakutani \cite{Kakutani, Kakutani72} has determined the diffraction for the weighting
system $w=(1,0)$ and it is 
\[\frac{1}{4}\delta_0 + \mbox{sc} \, ,\]
where sc is a non-trivial singular continuous measure on $\ZZ$. 
On the other hand, with the weighting $w=(1,1)$ the elements of $X_{\ZZ}$
are all just the sequence $\ZZ$ as far as the autocorrelation is concerned, and the 
diffraction is $\delta_{\ZZ}$. From these it follows that the diffraction for a general 
weighting system is

\[
(\frac{w_a+w_b}{2})^2\delta_0+(\frac{w_a-w_b}{2})^2 \mbox{sc} \, .
\]
In view of our remarks in \S 6.1, the diffraction for
$X_{\RR}$ is
\[
(\frac{w_a+w_b}{2})^2\delta_{\ZZ}+(\frac{w_a-w_b}{2})^2 \mbox{scp}
\]
where scp is the periodization of $\RR$ of the measure sc on
$\TT$.

The dynamical system is also mixed, pure point plus
singular continuous \cite{Keane}. There is an obvious continuous involution
$\sim$ on $X_\ZZ$ that interchanges the $a$ and $b$ symbols.
$L^2(X_{\ZZ},\mu_{\ZZ})$ splits into the $\pm 1$-eigenspaces for $\sim$:
$L^2(X_{\ZZ},\mu_\ZZ)=L^2_+(X_{\ZZ})\bigoplus L^2_-(X_{\ZZ})$.
$L^2_+(X_{\ZZ})$ is the pure point part of $L^2(X_{\ZZ},\mu)$ and its
eigenvalues are all the numbers of the form $k/{2^n}$,
$ n=0,1,\dots ; 0\leq k<2^n$ (literally $\exp{(2 \pi i k/{2^n})}$). 
On the other hand $L^2_-(X_{\ZZ})$ is singular continuous.

When we move to the suspension of $X_\ZZ$
we obtain $L^2(X_{\RR},\mu)$ which we know certainly retains
the eigenvalues of $L^2(X_{\ZZ},\mu_{\ZZ})$ and, due to the embedding
of $L^2(\RR,\widehat{\palm^{(1,0)}_1})$, also retains a singular continuous component. 

The dynamical spectrum is, of course,  independent of any particular assignments
of weights to a and b. We can draw the following conclusions from this:

{(i)} $w_a=1,w_b=0$.
\[\widehat{\palm^{(1,0)}_1}=\frac{1}{4}\delta_{\ZZ}+ \mbox{scp} \,.\]

The eigenfunctions of $L^2(\RR,\widehat{\palm^{(1,0)}_1})$ are
$\chF_{\{k\}},k\in\ZZ$. It follows that $\theta^w(\chF_{\{-k\}})$ is an
eigenfunction for eigenvalue $k$ (Thm.~\ref{main}, Cor.~ \ref{acscpp}). Thus
$\theta^w$ covers only the eigenvalues $k\in\ZZ$ of
$L^2(X_{\RR})$ and none of the fractional ones
$k/{2^n},n>0$. This shows that $\theta^w$ is not
surjective. Also $\theta^w$ embeds the
singular continuous part of $L^2(\RR,\widehat{\palm^{(1,0)}_1})$ into
$L^2_-(\RR)$, although we do not know the image.

{(ii)} $w_a=w_b=1$. In this case the diffraction is
$\delta_{\ZZ}$ (the Thue-Morse sequence with equal weights
looks like $\ZZ$). Although $L^2(\RR,\widehat{\palm^{(1,1)}_1})$ is pure point, its
image does not cover the pure point part of $L^2(X_{\RR},\mu)$, nor does it even generate it as an algebra. 
This shows that the requirement of unequal weights in Theorem~\ref{AlgebraGenerates} is
necessary.

{(iii)} $w_a=1,w_b=-1$. This time the diffraction is singular
continuous and 
$\theta^w$ does not even cover anything
of the the pure point part of $L^2(X_{\RR},\mu)$.

Cases (ii) and (iii) show that the non-existence of a particular
component in the diffraction spectrum implies nothing about its
existence or non-existence in the dynamical spectrum.

\subsection{Rudin-Shapiro}

We define the Rudin-Shapiro sequence using the notation of  ~\cite{Priebe}.
Consider the substitution rule $s$ defined on the alphabet $A'
:=\{1,\bar{1},2,\bar{2}\}$ as follows: $s(1)=1\bar{2},
s(2)=\bar{1}\bar{2}, s(\bar{1})=\bar{1}2, s(\bar{2})=12$. Let $\xi$ be the $s$-invariant sequence that starts with the symbol $1$. We can reduce this to a $2$-symbol sequence $\xi'$ with alphabet $\{a,b\}$ by replacing the symbols with no over-bar 
by the letter $a$ and the others by the letter $b$. This $2$-symbol sequence is 
usually called the Rudin-Shapiro sequence \cite{Quef}, though Priebe-Frank uses this appellation for the  original $4$-symbol sequence.

Let us start with the $2$-symbol sequence, which results in the $2$-coloured minimal and ergodic dynamical hull  $(X_\ZZ(\xi'), \ZZ)$, as developed above. There is a natural involution on the dynamical system that interchanges $a$ and $b$. Once again we introduce a system of weights $ w = (w_a, w_b)$. 

Under the system of weights $(1,-1)$ it is well known that the diffraction measure of the elements of $X_\ZZ(\xi')$ is the normalized Haar measure on $\RR/\ZZ$ \cite{Quef}, Cor. VIII.5. Thus 
$L^2(\RR,\widehat{\palm^{(1,-1)}_1}) =L^2(\RR,\Vol)$, where $\Vol$ is Lebesgue measure on $\RR$. 

On the other hand, the weighting system $(1,1)$ reduces the elements of $X_\ZZ(\xi')$
to copies of the sequence $\ZZ$. So, just as in the case of the Thue-Morse sequence,
we can deduce the general formula for the diffraction:
\[
(\frac{w_a+w_b}{2})^2\delta_{\ZZ}+(\frac{w_a-w_b}{2})^2 \Vol \, .
\]

The spectral decomposition
of $L^2(X_\ZZ(\xi'))$ is of the form
\[ L^2(X_\ZZ(\xi')) \simeq H  \oplus Z(f) \]
where $H$ is the pure point part with one simple eigenvalue 
$\exp(2\pi i q)$ for each dyadic rational number $q = a/2^n$, where $a\in \ZZ$, $n = 0,1, 2,  \dots$ \cite{Dekking, JMartin}; and
$Z(f)$ is a cyclic subspace which is  equivalent to $L^2(\RR,\Vol)$. In other words, the dynamical spectrum is mixed with a pure point and an absolutely continuous part \footnote{Explicitly $f$ is the function on 
$X_\ZZ(\xi')$ which is defined by 
$f(\zeta) = 1\, \mbox{or} \, -1$ according as $\zeta(0)$ is $a$ or $b$. This can be
deduced from the main theorem of \cite{Priebe}, where the equivalent result for dynamical system arising from the $4$ symbol sequence gives two copies of  $L^2(\RR,\Vol)$), and then by dropping to the factor.}.
In any case, we see that $L^2(X_\RR(\xi'), \mu)$ contains a pure point part whose eigenvalues include all the dyadic rationals, and also an absolutely continuous part into which the absolutely continuous part of $L^2(\RR, \widehat{\palm^w_1})$
must map by $\theta^w$. 

The analysis now proceeds exactly as in the case of the Thue-Morse sequence, with the same three types of possibilities except now the 
singular continuous parts are replaced by absolutely continuous parts.  

\subsection{Regular model sets}
In this example we see that even when everything is pure point and there is
only one colour, still $\theta$ need not be surjective. 

Let $(\Rd, \Rd, L)$ be a cut and project scheme with projection mappings
$\pi_i$, $i=1,2$. Thus $L$ is a lattice in $\Rd \times \Rd$, the projection $\pi_1$
to the first factor is one-one on $L$, and the projection $\pi_2(L)$ of $L$ has dense image
in the second factor. Let $W$ be a non-empty compact subset which is the closure
of its own interior and a subset of the second factor. We assume that
the boundary of $W$ has Lebesgue  measure $0$. The corresponding model
set is 
\[\Lambda(W) = \{ \pi_1(t) \,:\, t\in L, \pi_2(t) \in W \} \,. \]
It is a subset of $\cDr$ for some $r>0$ and it is pure point diffractive \cite{Hof3, Martin, BM}.
The orbit closure $X =\overline{ \Rd + \Lambda(W)}$ is uniquely ergodic. Its autocorrelation $\gamma$, and hence its diffraction $\widehat \gamma$, is
the same for all $\Gamma \in X$. Furthermore, the diffraction is explicitly known:
\[\widehat{\palm_1} =  \widehat\gamma = \sum_{k\in L^0} a_ k\delta_{\pi_1(k)} \]
where $L^0$ is the $\ZZ$-dual lattice of $L$ with respect to the 
standard inner product on $\Rd\times \Rd \simeq \RR^{2d}$ and
\[a_k = \left|\widehat{\chF_W}(-\pi_2(k))\right|^2  \, .\]
For more on this see \cite{Hof3}. The main point is that $\widehat{\palm_1}(\pi_1(k)) = 0$ if and only if $a_k = 0$.

Likewise $L^2(X, \mu)$ is known and it is isometric in a
totally natural way by an $\Rd$-map  to $L^2(\RR^{2d} /\ZZ^{2d}, \nu)$,
 where $\nu$ is Haar measure on the torus. Thus the spectrum of $X$ is pure point
and the eigenvalues are precisely all the points of $L^0$. Thus
the mapping $\theta$ embedding the diffraction into the dynamics will
be surjective if and only if for all $k\in L^0$, $a_k \ne 0$.

Now it is easy to see that we can find model sets for our given cut and project scheme
for which fail to be surjective at any $k\in L^0$ that we wish, as long as $k\ne0$.
To do this take $W$ to be something simple like a ball centered on $0$ and for
each scaling factor $s>0$ let $\gL^{(s)} := \gL(sW)$. The intensities of the Bragg peaks become
\[a_k^{(s)} :=  s^2\left|\widehat{\chF_W}(-s\pi_2(k))\right|^2 \, .\]
Since $\widehat{\chF_W}$ is continuous and takes positive and negative
values on every ray through $0$ in $\Rd$, but altogether takes the value
$0$ only on a meagre set, we see that by choosing $s$ suitably we can arrange
either that $a_k^{(s)}$ vanishes at any preassigned non-zero $k \in L^0$ 
(and $\theta$ is not surjective) or that alternatively 
 $a_k^{(s)}$ vanishes nowhere on $L^0$ (and $\theta$ is a bijection).

\subsection{The necessity of non-zero weights in Thm.~\ref{AlgebraGenerates}.}

Let $\gL =(\gL_a,\gL_b)$ where 
\[ \gL_a = \{ z \in \ZZ\,:\, z \equiv 0 \; \mbox{or} \;2\mod 4\}, \quad
 \gL_b = \{ z \in \ZZ\,:\, z \equiv 3\mod 4\}\, .\]
Then $\gL$ is periodic with period $4$ and its hull -- that is, the closure of its $\RR$ translation orbit -- is $X \simeq \RR/4\ZZ$ (a conjugacy of dynamical systems
with the standard action of $\RR$ on $\RR/4\ZZ$). 
Thus $L^2(X,\mu)$, where $\mu$ is Haar measure on $\RR/4\ZZ$, has pure point spectrum with eigenvalues $\frac{1}{4}\ZZ$.

Let $(w_a,w_b)$ be a weighting system for $\gL$.
The autocorrelation is everywhere the same and is easily seen to be
\[ \palm_1 = \frac{1}{2}w_a^2 \delta_{2\ZZ} + \frac{1}{4}w_aw_b \delta_{1+4\ZZ}+
 \frac{1}{4}w_aw_b \delta_{-1+4\ZZ} +  \frac{1}{4}w_b^2 \delta_{4\ZZ} \, . \]
 
 The Fourier transform, that is the diffraction, is then given by
 \begin{eqnarray*}  \widehat{\palm_1} &=& \frac{1}{4}w_a^2 \delta_{\frac{1}{2}\ZZ} + 
 \frac{1}{16}w_aw_b \exp(-2 \pi i (\centerdot)) \delta_{\frac{1}{4}\ZZ}\\
 &+& \frac{1}{16}w_aw_b \exp(2 \pi i (\centerdot)) \delta_{\frac{1}{4}\ZZ}+
 \frac{1}{16}w_b^2 \delta_{\frac{1}{4}\ZZ}\\
 &=& \frac{1}{4}\{(w_a^2 +\frac{1}{2}w_aw_b +\frac{1}{4}w_b^2)\delta_\ZZ
 + (w_a^2 -\frac{1}{2}w_aw_b +\frac{1}{4}w_b^2)\delta_{\frac{1}{2} +\ZZ}\\
 &+&\frac{1}{4}w_b^2\delta_{\frac{1}{4} + \ZZ} 
 +\frac{1}{4}w_b^2 \delta_{-\frac{1}{4} + \ZZ}\} \, . 
\end{eqnarray*}

Now it is clear that the image of $\theta$ can only generate eigenfunctions
for the eigenvalues $\pm \frac{1}{4} +\ZZ$ if $w_b \ne 0$ (and then in fact it
does so, independent of the value of $w_a$).

\section{The square-mean Bombieri-Taylor conjecture} \label{SecBombieri-Taylor}

\begin{theorem} \label{Bombieri-Taylor}
(The square mean Bombieri-Taylor conjecture)
Let $(X,\Rd,\mu)$ be a uniformly discrete, multi-coloured stationary ergodic point process, and assume that $w$ is a system of weights. Then the following are equivalent \footnote{Limits here are taken in the
$L^2$-norm on $(X,\Rd,\mu)$. Recently D.~Lenz  \cite{Lenz} has given a point-wise version of the result.}
:
\begin{itemize}
\item[(i)]
\[\frac{1}{\Vol (C_R)} \sum_{x\in C_R} \gl^w(\{x\}) \, e^{2 \pi i k.x} \nrightarrow 0
\quad \mbox{\rm as} \;  R\to\infty \; ;\]
 
\item[(ii)] $\widehat{\palm^w_1}(\{k\}) \ne 0$;
\item[(iii)] $k$ is an eigenvalue of $U$. 
\end{itemize}
In the case that $k$ is an eigenvalue, then
\[\frac{1}{\Vol (C_R)} \sum_{x\in C_R} \gl^w(\{x\}) \, e^{2 \pi i k.x} \to \theta^w(\chF_{k})\, .\]

\end{theorem}
For notational simplicity we shall prove the two technical lemmas that precede
the main proof in the $1$-dimensional case. However, it is easy to generalize
the proof to any dimension $d$. Throughout, $R$ is assumed to be a positive
integer variable.

\begin{lemma}
For all $\epsilon >0$,
\[ \lim_{\epsilon \to 0} \widehat{\palm^w_1}(\chF_{[-\epsilon,\epsilon]})=
\widehat{\palm^w_1}(\{0 \}) \, ,
\]
i.e. $\{\chF_{[-\epsilon,\epsilon]}\}_{\epsilon\searrow 0} \longrightarrow
\chF_{\{0\}}$ in $L^2(\RR, \widehat{\palm^w_1})$.
\end{lemma}

{\sc Proof:} Assume $\epsilon \to 0^+$. Let $F_{\epsilon} :=
\chF_{[-\epsilon,\epsilon]} - \chF_{\{0\}}$. Then for all $x\in
\RR$, $0\le F_\epsilon(x) \le 1$ and  $F_\epsilon(x)\searrow 0$
pointwise. Since $\widehat{\palm^w_1}$ is a translation bounded positive 
measure,  $\widehat{\palm^w_1}(F_\epsilon) \searrow 0$. Now,
\[ \int |  \chF_{[-\epsilon,\epsilon]} - \chF_{\{0\}}|^2 \mathrm{d}\widehat{\palm^w_1}
= \int F_\epsilon ^2 \mathrm{d}\widehat{\palm^w_1} \le  \int F_\epsilon
\mathrm{d}\widehat{\palm^w_1}  \longrightarrow 0 \, .\]  \qed

\begin{lemma} As functions of $y \in \RR$,
\[  \frac{1}{2R}\int_{-R}^R e^{2 \pi i y.x} \mathrm{d} x
\longrightarrow\chF_{\{0\}}(y)\]
in $L^2(\RR, \widehat{\palm^w_1})$ as $R \to \infty$.
\end{lemma}

{\sc Proof:} Let \[g_R(y):=  \frac{1}{2R}\int_{-R}^R e^{2 \pi i y.x} \mathrm{d} x =
\frac{\sin(2 \pi yR)}{2 \pi y R}\, .\]
We need to show that $\int_{-\infty}^\infty |g_R(y) - \chF_{\{0\}}(y) |^2
\mathrm{d}
\widehat{\palm^w_1}(y) \longrightarrow 0$. Since $|g_R(y) - \chF_{\{0\}}(y) | \le
F_\epsilon(y)$
for $-\epsilon \le y \le \epsilon$, we have
 $\int_{-\epsilon}^\epsilon |g_R(y) - \chF_{\{0\}}(y) |^2 \mathrm{d}
\widehat{\palm^w_1}(y) \longrightarrow 0$ as $\epsilon \to 0$, and the convergence is
uniform without reference to $R$.

For the remaining parts of the integral, we have (the part from $-\infty$ to
$-\epsilon$
is the same)
\begin{eqnarray*} \int_{\epsilon}^\infty |g_R(y) &-& \chF_{\{0\}} (y)|^2 \mathrm{d}
\widehat{\palm^w_1}(y) = \int_{\epsilon}^\infty \frac{\sin^2(2 \pi yR)}{(2 \pi y R)^2}
\mathrm{d} \widehat{\palm^w_1}(y)\\
&\le& \frac{1}{(2\pi R \epsilon)^2}\int_\epsilon^{\infty} \frac{\epsilon^2}{y^2}
\mathrm{d} \widehat{\palm^w_1}(y )\\
&\le& \frac{1}{(2\pi R \epsilon)^2}\left\{\int_\epsilon^{\epsilon+1} \mathrm{d}
\widehat{\palm^w_1}(y )
+\sum_{m=1}^\infty \frac{\epsilon^2}{m^2}\int_{\epsilon+m}^{\epsilon+m+1}
\mathrm{d}\widehat{\palm^w_1}(y )\right\} \, .
\end{eqnarray*}
Since $\int_{a}^{a+1}
\mathrm{d}\widehat{\palm^w_1}(y )$ is uniformly bounded by some constant $C(1)$ (due to the translation boundedness of $\widehat{\palm^w_1}$, 
Thm.~\ref{colourNautoCorr})
we see that 
\[ \int_{\epsilon}^\infty |g_R(y) - \chF_{\{0\}}(y) |^2 \mathrm{d}
\widehat{\palm^w_1}(y) \longrightarrow 0\] as long as $R\epsilon \to \infty$ as
$R\to\infty$.
Putting $\epsilon = R^{-1/2}$ gives the necessary convergence of both parts.
\qed \smallskip

\noindent {\sc Proof Theorem}\;\ref{Bombieri-Taylor}: 
(ii) $\Leftrightarrow$ (iii): $k$ is an eigenvalue if and only if
$-k$ is an eigenvalue, $\widehat{\palm^w_1}(\{k\}) = \widehat{\palm^w_1}(\{-k\})$
for all $k$, and $k$ is an eigenvalue if and only if $\widehat{\palm^w_1}(\{k\})\ne 0$.

\noindent
(iii) $\Leftrightarrow$ (i):
Let $f_R:= \frac{1}{2R} \chi_k \chF_{[-R,R]}$ (see \eqref{chars}).
Then
\begin{eqnarray*}
\widehat{f_R} (y) &=& \frac{1}{2R}\int_{\RR} e^{-2 \pi i y.x} \chi_k(x)
\chF_{[-R,R]}(x) {\mathrm d} x\\
&=&\frac{1}{2R}\int_{-R}^R  e^{2 \pi i (k-y).x} {\mathrm d} x \\
[1mm]\\
&\longrightarrow& \chF_{\{0\}}(k-y) = \chF_{\{k\}} (y) \, ,
\end{eqnarray*}
the convergence being as functions of $y$ in $L^2(\RR, \widehat{\palm^w_1})$ as $R \to
\infty$.

Let $\phi_{-k} = \theta^w(\chF_{\{k\}})$.   Thus $\widehat{f_R} \to \chF_{\{k\}}$ implies that $\theta^w(\widehat{f_R}) \to \phi_{-k}$ in $L^2(X,\mu)$, so 
\[ \int_X |N^w_{f_R}(\gl) - \phi_{-k}(\gl)|^2 \dmu(\gl) \to 0 \, ,\] which  from
(\ref{countingFunctions}) gives
\[ \int_X \left| \frac{1}{2R} \sum_{x \in
[-R, R]} \gl^w(\{x\}) \, e^{2 \pi i k.x} - \phi_{-k}(\gl) \right|^2
\dmu(\gl) \to 0 \, . \]
Thus $\frac{1}{2R} \sum_{x \in[-R, R]} \gl^w(\{x\}) \, e^{2 \pi i k.x}$ converges in
square mean to $\phi_{-k}$. Furthermore by Cor.~\ref{spectrum}, $\phi_{-k}$
is a $\chi_{-k}$-eigenfunction for $T_t$ if $\widehat{\palm^w_1}(k) \ne 0$
and is $0$ otherwise.

If $\phi_{-k} =0$ then 
$\frac{1}{2R} \sum_{x \in[-R, R]} \gl^w(\{x\}) \, e^{2 \pi i k.x} =0$ \, $\mu$-a.e.
If $\phi_{-k} \ne 0$ then $\{\gl\,:\, \phi_{-k}(\gl) = 0 \}$ is a measurable
$T$-invariant subset of $\mu$, since $\phi_{-k}$ is an eigenfuction, so by the
ergodicity it is of measure $0$ or $1$. It must be the former. Now using
the Fischer-Riesz theorem \cite{Dieudonne}, there is a subsequence of
$\{\frac{1}{2R} \sum_{x \in [-R, R]} \gl^w(\{x\}) \, e^{2 \pi i k.x} \}_R$ which
converges pointwise $\mu$-a.e. to $\phi_{-k}$. Since $\phi_{-k}$ is almost everywhere not zero,  $\frac{1}{2R} \sum_{x \in[-R, R]} \gl^w(\{x\}) \, e^{2 \pi i k.x}
\nrightarrow 0$.\qed

\section{A strange inequality} \label{inequality}

Let $(X, \Rd, \mu)$ be a uniformly discrete stationary ergodic
point process (no colour). Assume that the point sets of $X$
have finite local complexity, $\mu$- a.s. This implies 
that the autocorrelation measure $\palm_1$ is 
supported on a closed discrete subset of $\gL -\gL$ for
any $\gL$ whose autocorrelation is $\palm_1$. Thus for $\gL \in X$ we have, $\mu$-almost surely,
\[\palm_1(t) = \lim_{R\to\infty} \frac{1}{\ell(C_R)} \card((-t + \gL)\cap \gL \cap C_R) \,.\]

\begin{prop} \label{strange}
For all $k, t \in \Rd$, 
\[\left |e^{2\pi i k\cdot t}-1\right
|\widehat{\palm_1}^{1/2}(k)\le 2(\palm_1(0)-\palm_1(t))\, .\]
\end{prop}

{\sc \bf Proof:} Let $k \in \Rd$. Then

\[ \frac{1}{\ell(C_R)}\sum_{x\in \gL\cap C_R}e^{-2\pi i k\cdot
x} \longrightarrow g_{k}\]  
in the norm of $L^{2}(X,\mu)$, where $g_k$ is an eigenfunction of $T$ for the
eigenvalue $k$ if $\widehat{\palm_1}(k) \ne 0$ and $0$ otherwise (Thm.~\ref{Bombieri-Taylor}).

Suppose $\widehat{\palm_1}(k)\neq 0$.  Let $t\in \Rd$. Since
\[
(T_{t}g_{k})(\gL)=g_{k}(-t+\gL)=
\lim_{R\rightarrow\infty}\frac{1}{\ell(C_R)}\sum_{x\in(-t+\gL)\cap
C_R}e^{-2\pi i k\cdot x}
\]
for almost all $\gL \in X$,
\begin{eqnarray*}
((T_{t}-1)g_{k})(\gL)
&=&\lim_{R\rightarrow\infty}\frac{1}{\ell(C_R)}\left\{ \sum_{x\in
(-t+\gL)\cap C_R}e^{-2\pi i k\cdot x}-\sum_{x\in \gL \cap
C_R}e^{-2\pi i k\cdot x}\right\}\\
&=&\lim_{R\rightarrow\infty}h_R(\gL) \, , \end{eqnarray*}
where
\[
h_R(\gL):=\frac{1}{\ell(C_R)}\left\{\sum_{x\in (-t+\gL) \backslash \gL\cap
C_R}e^{-2\pi i k\cdot x}-\sum_{x\in \gL \backslash
(-t+\gL)\cap C_R}e^{-2\pi i k\cdot x}\right\}.
\]
Thus
$h_R\rightarrow (e^{2\pi i k\cdot t}-1)g_{k}$ in the
$L^{2}$-norm on $X$.

Furthermore,
\begin{eqnarray*}
|h_R(\gL)|
&\le & \frac{1}{\ell(C_R)}\left(\sum_{x\in (-t+\gL) \backslash \gL\cap
C_R}\left |e^{-2\pi i k\cdot x}\right |+\sum_{x\in \gL
\backslash (-t+\gL)\cap C_R}\left |e^{-2\pi i k\cdot x}\right |\right)\\
&\le&\frac{1}{\ell(C_R)}\left(\sum_{x\in (-t+\gL) \backslash \gL\cap
C_R}1+\sum_{x\in \gL\backslash (-t+\gL)\cap C_R}1\right)\\
&=&\frac{1}{\ell(C_R)}\sum_{x\in (-t+\gL) \bigtriangleup \gL\cap
C_R}1 =2(\palm_1(0)-\palm_1(t)).
\end{eqnarray*}
Note that $\palm_1(0)\geq\palm_1(t)$.

With these preliminaries out of the way, the rest of the proof
is straightforward. Since $\mu$ is a finite measure,
 $h_R\rightarrow (e^{2\pi i k\cdot t}-1)g_{k}$ in the $L^{1}$ norm also. Then
there is a subsequence $\{h_{R_i}\}$ of $\{h_R\}$ which converges to
$(e^{2\pi i k\cdot t}-1)g_{k}$ point-wise almost everywhere (\cite{Cohn}, Sec.~3.1).

Using the dominated convergence theorem ($|h_{R}(\gL)|\le 2\palm_1(0)$), we have
\begin{eqnarray*}
\lefteqn{\int_{X}\left |e^{2\pi i k\cdot t}-1\right|^{2}\left
|g_{k}(\gL)\right |^{2}d\mu(\gL)}\\
&=&\lim_{R_i \rightarrow
\infty}\int_{X}\left |h_{R_i}(\gL)\right |^{2}d\mu(\gL)
\le \int_{X}\left |2(\palm_1(0)-\palm_1(t))\right |^{2}d\mu(\gL)\\
&=&4\left |(\palm_1(0)-\palm_1(t))\right |^{2}.
\end{eqnarray*}
Meanwhile, from Thm.~\ref{main}
\begin{eqnarray*}
\int_{X}\left |e^{2\pi i k\cdot t}-1\right|^{2}\left
|g_{k}(\gL)\right |^{2}d\mu(\gL)&=&\left |e^{2\pi i
k\cdot
t}-1\right |^{2}\int_{\Rd}\chF_{k}^{2}d\widehat{\palm_1}\\
&=& \left |e^{2\pi i k\cdot
t}-1\right|^{2}\widehat{\palm_1}(k).
\end{eqnarray*}
So
\[
\left |e^{2\pi i k\cdot
t}-1\right|\widehat{\palm_1}^{\frac{1}{2}}(k)\le 2\left
|(\palm_1(0)-\palm_1(t))\right | =  
 2(\palm_1(0)-\palm_1(t))\, . \]\qed

\remark This result numerically links three interesting quantities. If
$\gL \in X$ has autocorrelation $\palm_1$ then
for the set $P(\epsilon)$ of $\epsilon$-statistical almost periods of $\gL$, i.e. $t$ for which $\palm_1(0) -\palm_1(t) <\epsilon$, 
the Bragg peaks $I(a)$ of intensity greater than $a >0$, i.e. $k$ for which $\widehat{\palm_1}(k) >a$, can occur only at points $k$ which are $2\epsilon/\sqrt a$-dual to $P(\epsilon)$, 
i.e. $k$ for which $\left |e^{2\pi i k\cdot t}-1\right| < 2\epsilon/\sqrt a$ for all
$t\in P(\epsilon)$. If this latter quantity is less than or equal to $1/2$
and either of $P(\epsilon)$ or $I(a)$ is relatively dense, then the other one is a Meyer set
\cite{M} Thm.~9.1. Furthermore, Bragg peaks can occur only on the $\ZZ$-dual of
the statistical periods ($t$ for which $\palm_1(t)=\palm_1(0)$), a fact that is of course 
very familiar in the case of crystals.\footnote{We are grateful to Nicolae Strungaru for this last observation.} We note that the inequality seems to be optimal. The maximum values
of $\left |e^{2\pi i k\cdot t}-1\right|$ and $\widehat{\palm_1}^{\frac{1}{2}}(k)$ are
$2$ and $\palm_1(0)$ respectively, whereas the minimum value of $\palm_1(t)$ is
$0$.

\section{Patterns and pattern frequencies}

Let $(X,\Rd,\mu)$ be a multi-colour uniformly discrete stationary ergodic point process. It is of interest to define the frequency of finite colour patterns in $X$. This is made difficult because from the built in vagueness of the topology of $X$ we know that we should not be looking for exact matches of some given colour pattern $F$ of $\cDr^{(m)}$, but rather close approximations to it. In addition there is the problem of how to anchor $F$, in order to specify it exactly as we move it around. This leads us to always assume that $F$ contains $0$, and then to define a  {\bf pattern} in $X$ as a pair $(F, V)$ where $F = \cup_{i=1}^m (F_i,i)$ is a finite
subset of $\cDr^{(m)}$ with $0 \in F^\downarrow := \cup F_i$ and $V$ is a bounded measurable neighbourhood of $0$ in $\Rd$. For a pattern
$(F, V)$ we then define the collection of elements of $X$ that contain it as
\[ X_{F,V} = X_{(F,V)} := \{ \gL \in X \,:\, F \subset V+\gL  \}  \, ,\]
and write ${\bf 1}_{F,V}$ for ${\bf 1}_{X_{F,V}}$.

Throughout one should keep in mind that $F$ and $\gL$ are multi-colour
sets, our conventions are that translations are by elements of $\Rd$ and are
always on the left, and the inclusions take colour into account.

For any bounded region $B$ define 
\[ L_{F,V}(\gL, B) := \card\{x \in\gL^\downarrow \,:\, F \subset -x + V +  \gL \, , \, x-V+F \subset B \} \,.\]

An initial idea for the frequency of the pattern $(F,V)$ in a set $\gL \in X$ might be:
\begin{equation} \label{freq}
\frq(\gL,F,V) = \lim_{R\to\infty} \frac{1}{\Vol(C_R)}L_{F,V}(\gL, C_R) \, .
\end{equation} 

This definition is very sensitive to the boundary of $V$, as one can see
from the simple example below. In general we do not know how to prove that this limit exists, even almost everywhere in $X$. However, we can prove that for $V$ open or $V$ closed, if the limit does exist then it is, almost surely, given by the Palm measure of $X_{F,V}$, and this, we know, does exist. Thus we are led to define:

The {\bf frequency} of the pattern $(F,V)$ in $X$ is $\palm(X_{F,V})$.

 The connection with Palm measures comes because (as is easy to see from the van Hove property of expanding cubes)  
\[ \lim_{R\to\infty} \frac{1}{\Vol(C_R)}L_{F,V}(\gL, C_R) = 
\lim_{R\to\infty}\frac{1}{\Vol(C_R)}\sum_{x\in \gL^\downarrow \cap C_R} {\bf 1}_{F,V}(-x + \gL) \, .\]
The latter is the average over $\gL$ (where the weighting system is all $1$s) of 
${\bf 1}_{F,V}$, if it exists. 

\begin{prop} \label{frequencies} Let $(F,V)$ be a pattern with $V$ an open set. Then 
\[ \lfrq(\gL,F,V) = \palm(X_{F,V}) \]
$\mu$-almost surely for $\gL \in X$, where $\lfrq$ means that the $\liminf$ is taken
in {\rm  (\ref{freq})}. Similarly, if 
$(F,V)$ is a pattern with $V$ a closed set, then 
\[ \ufrq(\gL,F,V) = \palm(X_{F,V}) \]
$\mu$-almost surely for $\gL \in X$, where $\ufrq$ means that the $\limsup$ is taken
in {\rm (\ref{freq})}.
 \end{prop}
 
 \begin{lemma} $X_{F,V}$ is open if $V$ is bounded and open and closed if $V$ is bounded and closed.
 \end{lemma}

{\sc \bf Proof:} Let $V$ be open and let $\gL \in X_{F,V}$. Then $F\subset
V+\gL$.  Since $V$ is open and $F$ is finite, there is an $\epsilon >0$ so that for
each $f\in F$, with $f = v+ x$, where $v\in V$, $x\in \gL$  (there may be choices, but fix one choice $x$ for each $f$), $f+ C_\epsilon \subset V+ x$. Choose $R>0$ so that $-V +F \subset C_R$. 
Let $\gL' \in U(\overline{C_R},C_\epsilon)[\gL]$ and let $f = v+x \in F$, as above. Since
$ \overline{C_R}\cap\gL \subset C_\epsilon+\gL' $ and $x \in  \overline{C_R}\cap\gL $,
$x= c+x' $ where $x' \in \gL'$, $c\in C_\epsilon$. Then $c+f \in V+ x$, so $f \in V+x'$.

Since $f \in F$ was arbitrary, $F \subset V+\gL'$ and $\gL' \in X_{F, V}$. Thus the open neighbourhood $U(\overline{C_R},C_\epsilon)[\gL]$ of $\gL$ lies in $X_{F, V}$.

The argument for $V$ closed is similar.
 \qed
\smallskip 

 {\sc \bf Proof of Prop.~\ref{frequencies} (sketch):} 
 Consider the case when $V$ is open. Then $X_{F,V}$ is open and the value of any 
 regular measure at  $X_{F,V}$ can be approximated as closely as desired by
a compact set $K\subset X_{F,V}$. For any such $K$ we can find a continuous 
  function $f$ with ${\bf 1}_K \le f \le {\bf 1}_{X_{F,V}}$.
Using Prop. \ref{colouraverageThm}, where all weights are assumed equal to $1$, 
 we obtain that $\palm(f)$ is almost surely the average of $f$ on $\gL$ and, from the definition of
 $f$, that for any $\epsilon >0$ and for large enough $R$,
 \[\palm(K)\le  \palm(f) \le \lfrq(\gL, F, V) \le \frac{1}{\Vol(C_R)}\sum_{x\in \gL \cap C_R} {\bf 1}_{F,V}(-x + \gL) + \epsilon \, .\]
 
 Integrating over $X$ and using the Campbell formula we have, independent of $R$,
 \[ \palm(K) \le \int_{X} \lfrq(\gL, F, V) \dmu \le \palm(X_{F,V}) \, . \]
 Now since we can make $\palm(K)$ as close as we wish to $\palm(X_{F,V})$ ,
 we obtain both
 \[\palm(X_{F,V}) \le \lfrq(\gL, F, V) \quad \mbox{and} \quad \int_{X} \lfrq(\gL, F, V) \dmu= \palm(X_{F,V})\]
 From this $\palm(X_{F,V}) = \lfrq(\gL, F, V)$, $\mu$-almost everywhere.  
 
 The result for $V$ closed is similar, this time approximating by open sets from
 above. \qed
 
{\sc \bf Example:} Consider the usual dynamical system based on $\ZZ$:
$X(\ZZ) \simeq \RR/\ZZ$. Let $F:= \{0, 1/4\}$ and $V:= (-1/4,1/4)$. For any
$\gL = t+ \ZZ$, we have 
\[\frac{1}{2n}\sum_{u\in (t+\ZZ) \cap [-n,n]}  {\bf 1}_{F,V} (-u + t + \ZZ) = 0\]
while
\[\frac{1}{2n}\sum_{u\in (t+\ZZ) \cap [-n,n]}  {\bf 1}_{F,\overline V} (-u + t + \ZZ) \thickapprox 1 \,.\]
In this case we have $X_{F,V} = \emptyset$, $X_{F,\overline{V}} = X$ and 
\[ 0 = \frq(\gL, F,V) = \palm(X_{F,V}) < \palm(X_{F,\overline{V}})
= \frq(\gL, F,\overline{V}) =1 \, . \]

\section{Final comments}
After Thm.~\ref{correlationsDetermineMu}, it is natural to ask whether or not for an
 $m$-coloured stationary ergodic uniformly discrete point process $(X, \mu)$ 
{\em all} the correlations of $\mu$ are necessary in order
to determine it. Is is possible that only a finite number of them will suffice? In \cite{Deng}
it is shown that given any $n \ge 2$ there are $1D$ examples based on multi-step Markov
processes for which only the $2,3, \dots, n$-point correlations are required to determine
$\mu$. 

In \cite{LM} the pure point case is studied and Thm.~\ref{correlationsDetermineMu} 
is used to relate the correlations to the extinctions (missing Bragg peaks) in the diffraction of $(X, \mu)$. For example, it is shown in the $1$-colour case that if there are no extinctions then the $2$ and $3$-point correlations determine the measure $\mu$. This seems to be the generic situation for regular model sets based on real internal spaces. 

\section{Acknowledgments}

The authors would like to thank Michael Baake and Daniel Lenz for their interest and their insightful comments on this work during its preparation. We also thank Haik Mashurian and Jeong-Yup Lee for their careful reading of the manuscript. Their contributions have substantially improved the paper.

\end{document}